\documentclass[12pt,a4paper,reqno]{amsart}

\usepackage[english]{babel}
\usepackage{amssymb,latexsym,amsfonts,amsthm,upref,amsmath}
\usepackage[foot]{amsaddr}
\usepackage[margin=1in]{geometry} 
\usepackage[dvipsnames,x11names]{xcolor}
 \definecolor{myblue}{HTML}{003399}
\usepackage{hyperref}
\hypersetup{colorlinks,citecolor=myblue,filecolor=black,linkcolor=myblue,urlcolor=myblue}
\usepackage{enumerate}  
\usepackage{tikz}
\usepackage{float}
\usepackage{cleveref}
\usepackage{mathtools}
\usepackage{cite}
\usepackage{filecontents}
\makeatletter
\newcommand{\leqnomode}{\tagsleft@true}
\newcommand{\reqnomode}{\tagsleft@false}

\makeatother
\newtheorem*{thm*}{Theorem}
\newtheorem*{lem*}{Lemma}
\newtheoremstyle{prim}{}{}{\normalfont}{}{\bfseries}{.}{ }{}
\newtheoremstyle{stil}{}{}{\slshape}{}{\bfseries}{.}{ }{}
\theoremstyle{stil}
\newtheorem{thm}{Theorem}[section]
\newtheoremstyle{defi}{}{}{}{}{\bfseries}{.}{ }{}
\theoremstyle{defi}
\newtheorem{defn}[thm]{Definition}
\theoremstyle{defi}
\newtheorem{rem}[thm]{Remark}
\theoremstyle{stil}
\newtheorem*{mthm*}{Main Theorem}
\newtheorem*{kor*}{Corollary}
\newtheorem{pro}[thm]{Proposition}
\theoremstyle{stil}
\newtheorem{lem}[thm]{Lemma}
\theoremstyle{stil}
\newtheorem{kor}[thm]{Corollary}
\theoremstyle{prim}

\newenvironment{prf}{\noindent \textit{Proof.}}{\null\hfill$\qed$\hskip
2mm\vskip 2mm}

\newcommand{\ar}{{\rm A}_{ h}(R)}

\newcommand{\vr}{V(R)}

\newcommand{\Yg}{  {\rm Y}_{\hspace{-1pt}h} (\mathfrak{gl}_{N})}

\newcommand{\aqpp}{\mathcal{A}_{q,p}(\widehat{\mathfrak{gl}}_2)}

\newcommand{\modd}{ \,{\rm mod\,\,}}

\newcommand{\Uc}{ \mathcal{U} }
\newcommand{\Yht}{  {\rm Y}_{\hspace{-1pt}h} ( \overline{R} )}
\newcommand{\Yhtc}{ \wtld{{\rm Y}}_{\hspace{-1pt}h} (  \overline{R}  )}

\newcommand{\Drp}{  {\rm D}_{ h} (R)}
\newcommand{\Dr}{  {\rm D} (\wvr{R})}

\newcommand{\R}{\wvr{R}}

\newcommand{\vac}{ \mathrm{\boldsymbol{1}}}

\newcommand{\z}{\mathfrak{z}}

\newcommand{\CC}{\mathbb{C}}

\newcommand{\ZZ}{\mathbb{Z}}

\newcommand{\Lc}{\mathcal{L}}

\newcommand{\Sc}{\mathcal{S}}

\newcommand{\Tc}{\mathcal{T}}

\newcommand{\wtld}{\widetilde}

\newcommand{\wvr}{\overline}

\newcommand{\ot}{\otimes}
\newcommand{\ts}{\hspace{1pt}}
\newcommand{\qdet}{ {\rm qdet}\hspace{1pt}}
\newcommand{\tr}{ {\rm tr}}

\newcommand{\ndo}{\mathop{\mathrm{End}}}
\newcommand{\om}{\mathop{\mathrm{Hom}}}

\newcommand{\rez}{\mathop{\mathrm{Res}}}

\newcommand{\iotaopjd}{\mathop{\iota_{z_1,z_2}}}

\newcommand{\iotaopdj}{\mathop{\iota_{z_2,z_1}}}
\newcommand{\iotaxy}{\mathop{\iota_{x,y}}}
\newcommand{\iotaxixj}{\mathop{\iota_{x_i,x_j}}}
\newcommand{\iotaxiyj}{\mathop{\iota_{x_i,y_j}}}

\newcommand{\iotaopjdhdva}{\mathop{\iota_{z_1,z_2}}}
\newcommand{\iotaopdjhdva}{\mathop{\iota_{z_2,z_1}}}

\newcommand{\fand}{\quad\text{and}\quad}
\newcommand{\Fand}{\qquad\text{and}\qquad}

\newcommand{\non}{\nonumber}
\newcommand{\beq}{\begin{equation}}
\newcommand{\eeq}{\end{equation}}
\newcommand{\ben}{\begin{equation*}}
\newcommand{\een}{\end{equation*}}

\makeatletter
\def\smalloverbrace#1{\mathop{\vbox{\m@th\ialign{##\crcr\noalign{\kern3\p@}%
  \tiny\downbracefill\crcr\noalign{\kern3\p@\nointerlineskip}%
  $\hfil\displaystyle{#1}\hfil$\crcr}}}\limits}
\makeatother

\makeatletter
\def\smallunderbrace#1{\mathop{\vtop{\m@th\ialign{##\crcr
   $\hfil\displaystyle{#1}\hfil$\crcr
   \noalign{\kern3\p@\nointerlineskip}%
   \tiny\upbracefill\crcr\noalign{\kern3\p@}}}}\limits}
\makeatother

\makeatletter
\def\author@andify{%
  \nxandlist {\unskip ,\penalty-1 \space\ignorespaces}%
    {\unskip {} \@@and~}%
    {\unskip \penalty-2 \space \@@and~}%
}
\makeatother

\begin{document}

\title{On two families of   quantum vertex algebras of FRT-type}

\author{Lucia Bagnoli}
\address[L. Bagnoli]{Instituto Nacional de Matem\'{a}tica Pura e Aplicada, Rio de Janeiro, RJ, Brazil}
\email{lucia.bagnoli@impa.br}

\author{Marijana Butorac}
\address[M. Butorac]{Faculty of Mathematics, University of Rijeka, Radmile Matej\v{c}i\'{c} 2, 51000 Rijeka, Croatia}
\email{mbutorac@math.uniri.hr}
 
\author{Slaven Ko\v{z}i\'{c}}
\address[S. Ko\v{z}i\'{c}]{Department of Mathematics, Faculty of Science, University of Zagreb,  Bijeni\v{c}ka cesta 30, 10000 Zagreb, Croatia}
\email{slaven.kozic@math.hr}


\begin{abstract}
We consider two new families of   quantum vertex algebras which are associated with the type $A$ trigonometric $R$-matrix and  elliptic $R$-matrix of the eight-vertex model. 
We show that their $\phi$-coordinated representation theory is governed by the so-called FRT-operator, $h$-adically restricted operator satisfying the FRT-relation,
and we demonstrate some applications of this result.
Finally, in the elliptic case, we investigate the properties of the quantum determinant associated with the corresponding  quantum vertex algebra.
\end{abstract}

\maketitle

\allowdisplaybreaks


\section{Introduction}\label{intro}
\setcounter{equation}{0}
\numberwithin{equation}{section}
 {\em Quantum vertex operator algebras} ({\em quantum VOAs}) were introduced by Etingof and Kazhdan \cite{EK}, motivated by {\em deformed chiral algebras} of 
E. Frenkel and Reshetikhin \cite{FR2}.
The first examples of quantum VOAs, the {\em quantum affine VOAs}   established in \cite{EK},  were associated with rational, trigonometric and elliptic $R$-matrices  of type $A$ and   featured a quantum version of   locality, the $\mathcal{S}$-{\em locality property}, which takes the form of the {\em quantum current commutation relation} of Reshetikhin and Semenov-Tian-Shansky \cite{RS}.
Afterward, the extensive development of the quantum vertex algebra theory was advanced by the  ideas   of Li \cite{LiG1,Li, Liphi}, Bakalov and Kac \cite{BK}, Li, Tan and Wang \cite{LiTW0,LiTW} and, later on, by De Sole, Gardini and Kac \cite{DGK} and Boyallian and Meinardi \cite{BM}. 
More recently, the well-known Frenkel--Jing problem \cite{FJ} of associating quantum vertex algebra theory with quantum affine algebras was settled 
by using Li's $\phi$-{\em coordinated module theory} \cite{Liphi}
in the series of papers by Jing, Kong, Li and Tan \cite{JKLT0,JKLT} and Kong \cite{Kong, Kong2}.

In the prior study by the first and  third author  \cite{BK1}, we considered   quantum vertex algebras such that  their $\mathcal{S}$-locality property is in the form of the famous {\em Faddeev--Reshetikhin--Takhtajan-relation}
({\em FRT-relation})  \cite{FRT}.
More specifically, we constructed a family of   such   quantum vertex algebras $\vr$, associated with the suitably normalized {\em Yang $R$-matrix} $R=R(u)$, and demonstrated that constructing a $\vr$-module 
over a topologically free $\CC[[h]]$-module $W$
 is equivalent to constructing a {\em FRT-operator over $W$}, i.e. a map $\mathcal{L}(u)\in \ndo\CC^N\ot \om(W,W((u))_h)$ such that 
\beq\label{frteqref}
R (u-v)\ts \mathcal{L}_1(u)\ts \mathcal{L}_2(v)=
 \mathcal{L}_2(v)\ts \mathcal{L}_1(u)\ts R (-v+u).
\eeq
This established  connections between   $\vr$ and the Yangian quantization of the Poisson algebra $\mathcal{O}(\mathfrak{gl}_N((z^{-1})))$,   which naturally occurs in the study of antidominantly shifted Yangians; see, e.g.,   the papers by  
 Frassek,  Pestun and Tsymbaliuk \cite{FPT} and Krylov and Rybnikov \cite{KR}.

In the present paper, we consider the type $A$ trigonometric $R$-matrix and  the elliptic $R$-matrix of the eight-vertex model, both   defined over the commutative ring $\CC[[h]]$ and denoted by $R(z)$; see Section \ref{sec01} for more details. 
In Section \ref{sec03}, we extend the aforementioned construction of  $\vr$ from \cite{BK1} to the trigonometric and elliptic case, so that the 
 $\Sc$-locality of these two new families of    quantum vertex algebras is again governed by the FRT-relation. 
Next, in Section \ref{FRTsection}, by adapting \eqref{frteqref} to the trigonometric and elliptic setting, we introduce the concept of {\em FRT-operator} as a map
$\Lc(z) \in\ndo\CC^N \ot \om(W, W ((z))_h)$ 
which satisfies
\beq\label{frteqref2}
   \iotaxy (R (x/y)) \ts \Lc_1 (x)\ts \Lc_2(y)
=
\Lc_2(y)\ts \Lc_1 (x)\ts  R^* (x/y). 
\eeq

In Section \ref{sec04_trig}, we restrict our considerations to the trigonometric case.
We show that constructing a FRT-operator over the topologically free $\CC[[h]]$-module $W$ is equivalent to constructing a structure of $\phi$-coordinated $\vr$-module over $W$ (with $\phi(z_2,z_0)=z_2 e^{z_0}$). Next, we demonstrate some applications of this result to the representation theory of {\em generalized $h$-Yangian extensions} and {\em $h$-Yangians}, certain classes of algebras based on \cite{BK1} and \cite{FPT}, respectively, which both extend the  ordinary $h$-Yangian associated with the trigonometric $R$-matrix of type $A$.
Moreover, we give a vertex operator-theoretic interpretation of certain commutative families and central elements of the completed $h$-Yangian found in \cite[Sect. 8]{BK1}.

In Section \ref{sec04_ell}, we consider the elliptic setting. Due to the specific form of the elliptic $R$-matrix, one needs to suitably adapt the notion of $\phi$-coordinated $\vr$-module   in order to obtain an analogue of the aforementioned result from the previous section. More specifically, we introduce the notion of 
{\em generalized} $\phi$-coordinated $\vr$-module
 (again with $\phi(z_2,z_0)=z_2 e^{z_0}$), which is based on the ideas of Li, Tan and Wang \cite{LiTW0,LiTW} and, in particular, on their notion of {\em ($\phi$-coordinated) twisted module}. Finally, we prove that   constructing a FRT-operator over the topologically free $\CC[[h]]$-module $W$  is equivalent to constructing a structure of generalized $\phi$-coordinated $\vr$-module over $W$. Although the proof of this result   relies on the usual $R$-matrix techniques and partially resembles the proof   \cite[Sect. 3]{K}, it also contains numerous technical arguments,  so we devote entire Subsection \ref{sec_prf} to presenting it in detail. 

In Section \ref{last_section}, we again consider the elliptic setting, where the identity \eqref{frteqref2}   can be regarded as a specialization of the defining relation
for the elliptic quantum algebra for $\widehat{\mathfrak{gl}}_2$ of
 Foda,   Iohara,   Jimbo,   Kedem,   Miwa and Yan \cite{FIJKMY}. Motivated by their notion of quantum determinant for the algebra $\mathcal{A}_{q,p}(\widehat{\mathfrak{gl}}_2)$, as well as by  its generalization to $\mathcal{A}_{q,p}(\widehat{\mathfrak{gl}}_N)$ by Frappat, Issing and Ragoucy \cite{FIR}, we define the {\em quantum determinant} for $\vr$. We show that its coefficients belong to the center of the  quantum vertex algebra $\vr$ and demonstrate its application to constructing elements of the submodule of invariants of a generalized $\phi$-coordinated $\vr$-module.

\section{Preliminaries}\label{sec01}

In this section, we recall the  trigonometric and     elliptic $R$-matrix of type $A$.  We use the same notation for both $R$-matrices,    which should not cause any confusion as
each   subsection  is devoted to   only  one of them.
In the rest of the paper, we   consider  both $R$-matrices over the ring $\CC[[h]]$ and use
the setting from Subsections \ref{sec0101} and   \ref{sec0103sec0103}. 
As for Subsections \ref{sec0101sec0101} and    \ref{sec0102sec0102}, they contain some intermediate definitions and results which  we needed to adapt the elliptic $R$-matrix 
 to the $h$-adic setting.

\subsection{Trigonometric \texorpdfstring{$R$}{R}-matrix over \texorpdfstring{$\CC[[h]]$}{C[[h]]}}\label{sec0101}

In this subsection, we recall the definition  of the  type $A$  trigonometric $R$-matrix. We consider   its additive
and multiplicative form, which both naturally occur in the theories of
quantum groups and exactly solvable models; see, e.g.,  \cite{FRT,J,PS}. In the setting of this paper, the former  governs the structure of a certain    quantum vertex algebra while the latter of its $\phi$-coordinated modules.

Let $N\geqslant 2$ be an integer and $h$ a formal parameter.
Consider  the two-parameter trigonometric $R$-matrix   $\R (z_1,z_2)\in\ndo\CC^N \ot \ndo\CC^N  [[h]][z_1,z_2]$ given by
\begin{align}
\R (z_1,z_2) =&\left(z_1 e^{-h/2}-z_2 e^{h/2}\right)\sum_{i=1}^N e_{ii}\ot e_{ii} 
+ (z_1-z_2)\sum_{\substack{i,j=1\\i\neq j}}^N e_{ii}\ot e_{jj} \non\\
&+ \left(e^{-h/2}-e^{h/2}\right)z_1 \sum_{\substack{i,j=1\\i> j}}^N e_{ij}\ot e_{ji}
+\left(e^{-h/2}-e^{h/2}\right)z_2\sum_{\substack{i,j=1\\i< j}}^N e_{ij}\ot e_{ji},\label{R2p}
\end{align}
where $e_{ij}$ are the matrix units. 
We shall also need its one-parameter counterparts, which are defined as follows.
First, introduce the power series
\begin{align}
&g(z)=  -e^{h/2}\left( 1-ze^{h}\right)^{-1}=-e^{h/2}\sum_{a\geqslant 0}z^a e^{ah}\in\CC[[z,h]],\label{g_trig}\\
&g(e^{u})=g(z)\big|_{z=e^u}= -e^{ h/2}\left(1-e^{u+h}\right)^{-1}\in\CC((u))[[ h]],\non
\end{align}
where we use the common expansion conventions, so that for any $k\in\CC$ we have
\begin{gather*}
e^{u+kh}=\sum_{a\geqslant 0}\frac{(u+kh)^a}{a!}\in\CC[[u,h]],\qquad
\frac{1}{u+kh}=u^{-1}\sum_{a\geqslant 0}\left(-\frac{kh}{u}\right)^a\in\CC[u^{-1}][[ h]],\\
\left(1-e^{u+kh}\right)^{-1}=\frac{1}{u+kh}\left(\frac{1-e^{u+kh}}{u+kh}\right)^{-1}
\in\CC((u))[[ h]] ,
\end{gather*}
with the inverse  of
$$
 \frac{1-e^{u+kh}}{u+kh} =-\sum_{a\geqslant 1}\frac{(u+kh)^{a-1}}{a!}\in \CC[[u,h]]
$$
being taken in $\CC[[u,h]]$.  Finally, introduce the normalized  $R$-matrices
\beq\label{R1p}
R(z)=g(z) \ts \R(z,1)
\fand
R(e^u)=g(e^u)\ts \R(e^u,1).
\eeq 

The $R$-matrix $R(z)$ (resp. $R(e^u)$)   belongs to
\beq\label{ptrs1}
\ndo\CC^N \ot  \ndo\CC^N  [[z,h]]\quad
\text{(resp. }\ndo\CC^N \ot  \ndo\CC^N   ((u))[[h]]\text{).}
\eeq
Both $R$-matrices satisfy the {\em quantum Yang--Baxter equation},  
\begin{align}
R_{12}(z_1)\ts R_{13}(z_1 z_2)\ts R_{23}(z_2)
&=   R_{23}(z_2)\ts R_{13}(z_1 z_2)\ts R_{12}(z_1),\label{ybe1}\\
R_{12}(e^{u_1})\ts R_{13}(e^{u_1+u_2})\ts R_{23}(e^{u_2})
&=  R_{23}(e^{u_2})\ts R_{13}(e^{u_1+u_2})\ts R_{12}(e^{u_1}).\label{ybe2}
\end{align}
Moreover, they possess the {\em unitarity properties}
\beq\label{uni}
R_{12}(z)\ts R_{21}(1/z)=1 \fand R_{12}(e^u)\ts R_{21}(e^{-u})=1.
\eeq

To consider the trigonometric and elliptic setting simultaneously, from now on  we write
$R^*(z)=R(z)$ and $R^*(e^u)=R(e^u)$
(in the elliptic case we have \eqref{zvijezdice}). 
As the aforementioned $R$-matrices   possess the unitarity property \eqref{uni}, we also set 
$
\Uc(z)=\Uc(e^u)=1
$
 (in the elliptic case we have \eqref{runi}).

 \begin{rem}\label{rationalfs_trig}
Recall the formal Taylor Theorem,
$$
a(z+z_0)=\sum_{k=0}^{\infty}\frac{z_0^k}{k!} \frac{\partial^k}{\partial z^k} a(z)\quad\text{for}\quad a(z)\in V[[z^{\pm 1}]],
$$
where $V$ is a vector space.   
Due to  the expansion
\beq\label{taylor}
\frac{1}{1-ze^{  u- v+ a  h}}=\sum_{k=0}^{\infty} \frac{\left(e^{  u-v+ a  h}-1\right)^k z^k}{k!}
\frac{\partial^k}{\partial z^k} \left(\frac{1}{1-z}\right), \quad\text{where }a \in\CC,
\eeq
and the form \eqref{g_trig} of the   series $g(z)$,
we can regard  the $R$-matrix $R(ze^{  u- v+ a  h})$   as a rational function  in  $z$, i.e. as an element  of
$
  (\ndo\CC^N)^{\ot 2}   (z)[[u,v,h ]] .$
\end{rem}

 It is evident  from the  form  of the $R$-matrix that it satisfies the following proposition.
	\begin{pro}\label{prop_trig_novi}
	\begin{enumerate}
\item	For any integer $k\geqslant 0$ there exists an integer $r\geqslant 0$ such that
$$
(1-z )^r R(z)^{\pm 1}\in
\ndo\CC^N  \ot \ndo\CC^N \ts [z  ,   h]
\mod h^k.
$$

\item
For any  integers $a,b,k > 0$ and $\alpha\in\CC$ there exists an integer $r\geqslant 0$ such that the coefficients of all monomials 
\beq\label{xmonomix_trig}
u^{a'} v^{b'} h^{k'} , \quad\text{where}\quad 0\leqslant a'\leqslant a,\, 0\leqslant b'\leqslant b,\, 0\leqslant k'\leqslant k ,
\eeq
 in
$ (z_1  - z_2 )^r R (z_1 e^{u -v +\alpha h} / z_2)^{\pm 1}$ 
belong to $(\ndo\CC^N)^{\ot 2}[z_1  ,z_1^{- 1 }, z_2 , z_2^{- 1 }]$ and such that the coefficients of all monomials \eqref{xmonomix_trig} in
$$
\left((z_1  - z_2 )^r R (z_1 e^{u -v +\alpha h} / z_2)^{\pm 1}\right)\Big|_{z_1 =z_2  e^{z_0 }}^{\modd u^{a}, v^{b}, h^{k} }\big. 
\fand   
z_2^{ r} (e^{ z_0}-1)^r  R (e^{z_0+u-v+\alpha h})^{\pm 1}
$$
coincide.
\end{enumerate}
\end{pro}

\subsection{Elliptic \texorpdfstring{$R$}{R}-matrix}\label{sec0102}

\subsubsection{\texorpdfstring{$R$}{R}-matrix over \texorpdfstring{$\mathbb{C}(q^{1/2})[[p^{1/2}]]$}{C(q1/2)[[p1/2]]}}\label{sec0101sec0101}

Let us recall the definition of the elliptic $R$-matrix of the 
 eight-vertex
  model. For more information on the elliptic $R$-matrices see  \cite{Bax,Bel,Chu, RT,T}. 
	We  use the formal power series point of view, i.e. we  regard the $R$-matrix parameters $z$, $p^{1/2}$, $q^{1/2}$ as indeterminates.
	Despite some notational differences, our exposition closely follows Foda,   Iohara,   Jimbo,   Kedem,   Miwa and Yan \cite[Sect. 2]{FIJKMY}. In particular, the $R$-matrix $R(z)$,   defined by \eqref{rbar} below, coincides with   $R^+(\zeta)$ from \cite[Sect. 2]{FIJKMY}.

 For a variable $z$ and a parameter  $p $ define
$$
(z;p )_\infty
=\prod_{k\geqslant 0} (1-zp^{k} ) \in \CC[[z,p]].
$$
It is known by \cite{FR} that there exists a unique   power series $f(z)\in\CC(q)[[z]]$ of the form
$$
f(z)=1+\sum_{k\geqslant 1} f_k z^k
\quad\text{such that}\quad
f(z q^{4}) =f(z)\frac{(1-zq^2)^2 }{(1-z)(1-zq^{4})}.
$$
It is given by
\beq\label{f}
f(z)=\frac{(z;q^{4})_\infty\ts (zq^{4};q^{4})_\infty}{(zq^{2};q^{4})_\infty^2}.
\eeq

Introduce the power series  $\alpha(z),\beta(z)\in\CC(q^{1/2})[z^{\pm 1}][[p^{1/2}]]$ by
\begin{align}\label{alphabeta}
\alpha(z) = \frac{(p^{1/2}qz;p)_\infty}{(p^{1/2}q^{-1}z;p)_\infty}\prod_{k\geqslant 1} f(p^k z^2)\fand
\beta(z)= \frac{(-pqz;p)_\infty}{(-pq^{-1}z;p)_\infty}
\prod_{k\geqslant 1} f(p^k z^2).
\end{align}
The elliptic $R$-matrix  is given by
\beq\label{rbar}
R(z)=\begin{pmatrix}
a(z)& & & d(z)\\
 & b(z) & c(z) & \\
 & c(z) & b(z) & \\
d(z)& & & a(z)
\end{pmatrix}.
\eeq
Its matrix entries are  
  elements of $\CC(q^{1/2})((z))[[p^{1/2}]]$ such that $a(z)$ and $b(z)$ possess only even and   $c(z)$ and $d(z)$ only  odd powers of $z$. They are  uniquely determined by  
\begin{align}\label{prghj3}
a(z)+d(z)=q^{-1/2}f(z^2)^{-1}\frac{\alpha(z^{-1})}{\alpha(z)}\fand
b(z)+c(z)=q^{1/2}f(z^2)^{-1} \frac{1+q^{-1}z}{1+qz}\frac{\beta(z^{-1})}{\beta(z)}.
\end{align}
The $R$-matrix \eqref{rbar} satisfies
\beq\label{runiii}
R_{12}(z)\ts R_{21}(z^{-1})=\Uc(z), \quad\text{where}\quad \Uc(z)=  
 q^{-1} f(z^2)^{-1} f(z^{-2})^{-1}. 
\eeq

\subsubsection{\texorpdfstring{$R$}{R}-matrix over \texorpdfstring{$\CC[[h,p^{1/2}]]$}{C[[h,p1/2]]}}\label{sec0102sec0102}
 
Our next goal is to show that, instead of the field $\CC(q^{1/2})$, the elliptic $R$-matrix $R(z)$ can be  regarded over the commutative ring $\CC[[h]]$, where $h$ is the indeterminate such that
\beq\label{qeh}
q^{1/2}=e^{h/4} =\sum_{k\geqslant 0} \frac{h^k}{4^k k!} \in\CC[[h]].
\eeq
Moreover, we shall consider its additive form over $\CC[[h]]$, which is obtained via
\beq\label{zeu}
z=e^u =  \sum_{k\geqslant 0} \frac{u^k}{k!}\in\CC[[u]] .
\eeq

As demonstrated in \cite[Sect. 2]{KM}, the series \eqref{f} admits the presentation
\beq\label{in4}
f(z)=1+\sum_{k\geqslant 1} a_k \left(\frac{z}{1-z}\right)^k\quad\text{for some}\quad a_k\in\CC(q)
\eeq
such that all $a_k (1-q)^{-k}\in\CC(q)$ are regular at $q=1$. Thus,   applying the substitution $q=e^{h/2} $ to \eqref{in4},    we obtain
\beq\label{in5}
f(z)\big|_{q=e^{h/2}}\big.=1+\sum_{k\geqslant 1} b_k \left(\frac{z}{1-z}\right)^k \quad\text{with}\quad b_k = a_k\big|_{q=e^{h/2}}\big.\in h^k\CC[[h]].
\eeq

Let $\CC_*(u)$ be the algebra extension of $\CC[[u]]$ obtained by inverting all  nonzero polynomials. Denote by $\iota_u$ the canonical embedding $\CC_*(u)\to\CC((u))$.
Applying the substitution $z=e^u $ to \eqref{in5}  and then taking the image under $\iota_u$  we obtain
\beq\label{in3}
 \iota_u \ts f(z)\big|_{q=e^{h/2},z=e^u}\big. \in\CC((u))[[h]].
\eeq
Actually,  in \eqref{in3}, we  applied an extension of   $\iota_u$,   the  embedding  $\CC_*(u)[[h]]\to \CC((u))[[h]]$. In what follows, we  often employ the suitable extensions and generalizations of the map $\iota_u$.  To simplify the notation,
we usually omit the symbol $\iota_u$. For example,
 we denote the series in \eqref{in3} by $f(e^u)$. 
It is worth noting that     \eqref{in3}   is invertible in $\CC((u))[[h]]$.

Let $c\in\CC$. Denote by  $R^*(z)$     the $R$-matrix obtained by replacing the parameter $p$   with $p^*\coloneqq pq^{-2c}$ in $R(z)$,  
\beq\label{zvijezdice}
R^*(z)=R(z)\big|_{p=p^*}.
\eeq

\begin{lem}\label{l1}
\begin{enumerate}
\item
By applying   the substitution \eqref{qeh}   to  the $R$-matrices  $R(z)$ and $R^*(z)$
we obtain   well-defined elements   
\beq\label{ermatrice2}
R(z)\big|_{q=e^{h/2} }
,\,
R^*(z)\big|_{q=e^{h/2} }\in
  \ndo\CC^2\ot \ndo\CC^2  [[h]]((z ))[[p^{1/2}]]. 
\eeq
\item
By applying the substitutions
\eqref{qeh} and \eqref{zeu}
and then the embedding $\iota_u$   to     $R(z)$ and $R^*(z)$
we obtain   well-defined elements 
\beq\label{ermatrice1}
\iota_u\ts R(z)\big|_{q=e^{h/2},z=e^u}
,\,
\iota_u\ts R^*(z)\big|_{q=e^{h/2},z=e^u}\in   \ndo\CC^2\ot \ndo\CC^2  [[p^{1/2}]]((u))[[h]].
\eeq
\end{enumerate}
\end{lem}

 \begin{prf}
The lemma   is verified by examining the explicit expression for the $R$-matrix. We only consider the $R(z)$ case, as the proof for $R^*(z)$ goes analogously. 
Let us prove \eqref{ermatrice2}.
First, due to the form of the power series   $\alpha(z)$ and $\beta(z)$, as given by  \eqref{alphabeta}, we have
\beq\label{prghj1}
q^{-1/2}\frac{\alpha(z^{-1})}{\alpha(z)},\, q^{1/2}  \left(1+q^{-1}z\right)\frac{\beta(z^{-1})}{\beta(z)}\in\CC[q^{\pm 1/2}][z^{\pm 1}][[p^{1/2}]] .
\eeq
Next, it is clear that
\beq\label{prghj2}
(1+qz)^{-1} \in\CC[q][[z]].
\eeq
Thus, since $f(z^2)^{-1}\in\CC(q)[[z]]$,   we conclude by \eqref{prghj1} and \eqref{prghj2} that the right-hand sides of both equalities in \eqref{prghj3} belong to $ \CC[q^{\pm 1/2}]((z ))[[p^{1/2}]] $. Hence     the $R$-matrix entries $a(z)$, $b(z)$, $c(z)$, $d(z)$, which appear on the left-hand sides of these equalities, are elements of $ \CC[q^{\pm 1/2}]((z ))[[p^{1/2}]] $ as well. Finally, by applying the substitution \eqref{qeh} to these entries, we find that the first assertion of the lemma holds for $R(z)$.  

As for the second assertion, by applying the substitutions \eqref{qeh} and \eqref{zeu} to the terms in \eqref{prghj1}, we obtain elements of
$\CC[[p^{1/2},u,h]]$. Therefore,  their image under the embedding  $\iota_u$ lies in $\CC[[p^{1/2},u,h]]$ as well. 
Regarding the term \eqref{prghj2}, we have 
$$
 \iota_u (1-q  z )^{-1}\big|_{q=e^{h/2},z=e^u}=
\iota_u (1-e^{u+h/2})^{-1}=(u+h/2)^{-1}\frac{u+h/2}{1-e^{u+h/2}} \in \CC((u))[[h]],
$$
  due to   the factor
	$(u+h/2)(1-e^{2u+h})^{-1}$ being invertible in $\CC[[u,h]]$ and the
	expansion
$$
(u+h/2)^{-1}=
\sum_{k\geqslant 0}\binom{-1}{k}  u ^{-1-k}   (h/2)^h   \in\CC[u^{-1}][[h]].
$$
Next, using \eqref{in5} and \eqref{in3}, one   checks that
$ 
 \iota_u \ts f(z^2)^{-1}\big|_{q=e^{h/2},z=e^u}\big. $ belongs to $ \CC((u))[[h]]
$.
The preceding discussion   implies that by applying the substitutions
\eqref{qeh} and \eqref{zeu}
and then the embedding $\iota_u$   to the right-hand sides of both equalities in \eqref{prghj3} we get elements of $\CC[[p^{1/2}]]((u))[[h]]$. The same holds for the right-hand sides of the equalities
\begin{align}\label{prghj6}
a(z)-d(z)=q^{-1/2}f(z^2)^{-1}\frac{\alpha(-z^{-1})}{\alpha(-z)}\fand
b(z)-c(z)=q^{1/2}f(z^2)^{-1} \frac{1-q^{-1}z}{1-qz}\frac{\beta(-z^{-1})}{\beta(-z)},
\end{align}
which are found by replacing $z$ with $-z$ in \eqref{prghj3}. Finally,  as the $R$-matrix entries  $a(z)$, $b(z)$, $c(z)$, $d(z)$ can be expressed as  linear combinations of the corresponding right-hand sides of the identities \eqref{prghj3} and \eqref{prghj6}, we conclude that the second assertion of the lemma holds for $R(z)$.
 \end{prf}

 \begin{rem}\label{rationalfs}
 By the proof of Lemma \ref{l1} and the expansion in \eqref{taylor},
we can regard  the $R$-matrices $R(ze^{  u-v+ a h})$ and $R^*(ze^{  u-v+ a h})$ with $a\in\CC$ as rational functions in  $z$, i.e. as  elements of
$
  (\ndo\CC^2)^{\ot 2}   (z)[[u,v,h,p^{1/2}]] .$
\end{rem}

\subsubsection{\texorpdfstring{$R$}{R}-matrix over \texorpdfstring{$\CC[[h]]$}{C[[h]]}}\label{sec0103sec0103}

The    $R$-matrix \eqref{rbar} (resp. \eqref{ermatrice2}) depends on the parameters $p,q$ (resp. $p,h$) and the variable $z$. Their  physical interpretation in terms of the  eight-vertex model can be found in the book by Baxter \cite[Chap. 10.4, 10.7]{Bax}. As  indicated by  \eqref{qeh}, throughout the rest of the paper, we set $q=e^{h/2}$, where $h$ is a formal parameter. In addition,  from now on, we set 
\beq\label{pandh}
p=ah^{2b}\quad\text{for}\quad a\in\mathbb{R},\, b\in\mathbb{Z},\, a,b>0.
\eeq
This makes the elliptic $R$-matrix compatible with the ($h$-adic) quantum vertex algebra theory \cite{EK,Li}, which  we want to associate with it.

Following the above convention,    we denote the $R$-matrices  in \eqref{ermatrice2} (resp. \eqref{ermatrice1})  by   $R(z)$ and $R^*(z)$  (resp. $R(e^u)$ and $R^*(e^u)$), respectively, where, from now on, we assume that the elliptic nome $p$ is replaced by $ah^{2b}$ as in \eqref{pandh}.
By Lemma \ref{l1} and the discussion following the Remark \ref{rationalfs},  
  $R(z)$ and $R^*(z)$  (resp. $R(e^u)$ and $R^*(e^u)$) belong to 
\beq\label{ptrs2}
\ndo\CC^2 \ot  \ndo\CC^2  ((z ))[[h]]\quad
\text{(resp. }\ndo\CC^2 \ot  \ndo\CC^2   ((u))[[h]]\text{).}
\eeq
The $R$-matrices $R(z)$ and   $R(e^u)$ satisfy the {\em quantum Yang--Baxter equation},  
\begin{align}
R_{12}(z_1)\ts R_{13}(z_1 z_2)\ts R_{23}(z_2)
&=   R_{23}(z_2)\ts R_{13}(z_1 z_2)\ts R_{12}(z_1),\label{rybe1}\\
R_{12}(e^{u_1})\ts R_{13}(e^{u_1+u_2})\ts R_{23}(e^{u_2})
&=  R_{23}(e^{u_2})\ts R_{13}(e^{u_1+u_2})\ts R_{12}(e^{u_1}).\label{rybe2}
\end{align}
Also, from now on,   $f(z)$    stands for  $f (z)\left|_{q=e^{h/2}}\right.\in\CC[[z,h]]$, as in \eqref{in5},  and we write 
\beq\label{useries2}
\Uc(z)=  
e^{-h/2} f(z^2)^{-1} f(z^{-2})^{-1}\fand 
\Uc(e^u)= \iota_u \,\ts\Uc(z)\big|_{q=e^{h/2},z=e^u};
\eeq
recall \eqref{runiii}.
Note that, due to \eqref{in3},      the  latter series is a well-defined element of $\CC((u))[[h]]$, which is
 even  with respect to the variable  $u$, so that we have $\ts\Uc(e^u)=\Uc(e^{-u})$.
 The $R$-matrices  $R(z)$ and  $R(e^u)$ possess the   property, which resembles unitarity, 
 \beq\label{runi}
R_{12}(z)\ts R_{21}(1/z)=\Uc(z)
\fand
R_{12}(e^u)\ts R_{21}(e^{-u})=\ts\Uc(e^u).
\eeq

We now extend the map $\iota$ to the case of multiple variables in fractional powers. Consider the ring  
	$\CC((h)) [[z^{1/2}_1,\ldots ,z^{1/2}_n]]$.
Denote 
by $\CC((h))_* (z^{1/2}_1,\ldots ,z^{1/2}_n)$
its localization    at nonzero polynomials $\CC((h)) [z^{1/2}_1,\ldots ,z^{1/2}_n]^{\times}$. There exists a   canonical embedding 
\beq\label{recall_eg}
\iota_{z_1,\ldots ,z_n}\colon	\CC((h))_* (z^{1/2}_1,\ldots ,z^{1/2}_n)\to \CC((h)) ((z^{1/2}_1))\ldots ((z^{1/2}_n)).
\eeq
We shall need several other   variations of this map, such as
\begin{align*} 
\CC (z_1,\ldots ,z_n)[[h]]&\to \CC ((z_1))\ldots ((z_n))[[h]],\\
\CC_* (z^{1/2}_1,\ldots ,z^{1/2}_n)[[h]]&\to \CC ((z^{1/2}_1))\ldots ((z^{1/2}_n))[[h]] 
\end{align*}
etc.,
 and we shall denote them by  $\iota_{z_1,\ldots ,z_n}$ as well. Whenever we apply such embeddings, Lemma \ref{l1} and the form of the elliptic $R$-matrix will ensure that their images possess only nonnegative powers of the parameter $h$.
In addition, we shall often omit the  symbol $\iota$ and employ  the usual expansion convention  where the embedding is determined by the order of the variables. For example, if $\tau$ is an element of the symmetric group $\mathfrak{S}_n$, then $f(e^{u_{\tau_1} + \ldots +u_{\tau_n}})$ denotes the series $\iota_{u_{\tau_1} , \ldots ,u_{\tau_n}}f(e^{u_{\tau_1} + \ldots +u_{\tau_n}})$ which belongs to  $\CC((u_{\tau_1}))\ldots ((u_{\tau_n}))[[h]]$.  Note that by this convention  we have, e.g.,  $f(e^{u_1+u_2})\neq f(e^{u_2+u_1})$,  as the former series belongs to $\CC((u_1))((u_2))[[h]]$ and the latter to $\CC((u_2))((u_1))[[h]]$.

By  following the   proof of Lemma \ref{l1}, one   checks that all its statements hold for the $R$-matrices $R(z^{-1})$ and   $R^*(z^{-1})$ as well. Hence, by   \eqref{runiii} and \eqref{useries2},   it follows that all assertions of Lemma \ref{l1}   hold for   $R(z)^{-1}$ and $R^*(z)^{-1}$ as well. 
The next proposition follows by an argument  which goes in parallel with the 
 proof of \cite[Lemma 3.2]{K} and relies on 
the forms of the series $f(z)$ and the $R$-matrix $R(z)$, along with Remark \ref{rationalfs}.
It can be regarded as an elliptic analogue of Proposition \ref{prop_trig_novi}.

\begin{pro}\label{localitycor}
	\begin{enumerate}
\item	For any integer $k\geqslant 0$ there exists an integer $r\geqslant 0$ such that
\beq\label{tlblm}
(1-z^2)^r R(z)^{\pm 1}\in
\ndo\CC^2  \ot \ndo\CC^2 \ts [z ,z^{- 1},  h]
\mod h^k.
\eeq

\item
For any  integers $a,b,k > 0$ and $\alpha\in\CC$ there exists an integer $r\geqslant 0$ such that the coefficients of all monomials 
\beq\label{xmonomix}
u^{a'} v^{b'} h^{k'} , \quad\text{where}\quad 0\leqslant a'\leqslant a,\, 0\leqslant b'\leqslant b,\, 0\leqslant k'\leqslant k ,
\eeq
 in
$ (z_1^2 - z_2^2)^r R (z_1 e^{u -v +\alpha h} / z_2)^{\pm 1}$ 
belong to $(\ndo\CC^2)^{\ot 2}[z_1 ,z_1^{- 1}, z_2, z_2^{- 1}]$ and such that the coefficients of all monomials \eqref{xmonomix} in
\beq\label{xmonomix222}
\left((z_1^2 - z_2^2)^r R (z_1 e^{u -v +\alpha h} / z_2)^{\pm 1}\right)\Big|_{z_1=z_2 e^{z_0}}^{\modd u^{a}, v^{b}, h^{k} }\big. 
\fand   
z_2^{2r} (e^{2z_0}-1)^r  R (e^{z_0+u-v+\alpha h})^{\pm 1}
\eeq
coincide.
\end{enumerate}
\end{pro}

Regarding the notation used in Proposition \ref{localitycor},
throughout the paper we often consider formal power series modulo sufficiently large powers of certain parameters or variables, which is  indicated by the ``mod'' symbol. 
For example, the assertion   \eqref{tlblm} means that if we write $(1-z^2)^r R(z)^{\pm 1}$ in the form 
$\sum_{s\geqslant 0} a_{s}^\pm(z)\ts h^s $ with $a_{s}^\pm(z) \in (\ndo\CC^2)^{\ot 2}   [[z^{ \pm 1 } ]]$, then  
$
\sum_{ k>s\geqslant 0 } a_{s}^\pm(z)\ts h^s $
belongs to 
$
\ndo\CC^2  \ot \ndo\CC^2 \ts [ z,z^{ - 1}, h]$.
Moreover, regarding the notation in   the second part of Proposition \ref{localitycor},
suppose we have  
$$
(z_1^2 - z_2^2)^r R (z_1 e^{u -v +\alpha h} / z_2)^{\pm 1}=\sum_{i,j,s,t\geqslant 0} a^{\pm }_{i,j,s}(z_1,z_2)\ts u^i v^j h^s  
$$
for some $a^{\pm }_{i,j,s}(z_1,z_2)\in (\ndo\CC^2)^{\ot 2} [[z_1^{\pm 1 } , z_2^{\pm 1 } ]]$.
The first assertion  implies that
$$
A^{\pm }_{a,b,k }\coloneqq\sum_{\substack{a\geqslant i\geqslant 0,\, b\geqslant j\geqslant 0\\k\geqslant s\geqslant 0 }} a^{\pm }_{i,j,s,t}(z_1,z_2)\ts u^i v^j h^s  
$$
belongs to $(\ndo\CC^2)^{\ot 2} [ z_1, z_1^{- 1 }, z_2, z_2^{- 1 },u,v,h ]$.
Hence  the substitution in \eqref{xmonomix222},
\beq\label{aswith}
\left((z_1^2 - z_2^2)^r R (z_1 e^{u -v +\alpha h} / z_2)^{\pm 1}\right)\Big|_{z_1 =z_2  e^{z_0 }}^{\modd u^{a}, v^{b}, h^{k}}\big. \coloneqq
A^{\pm }_{a,b,k}\Big|_{z_1 =z_2  e^{z_0 }},
\eeq
is well-defined.
In \eqref{aswith},   the ``mod'' symbol  indicates that the substitution  in the subscript  is not  applied to the entire expression inside the brackets, but rather to its coefficients with respect to the sufficiently small powers of the variables and parameters, as given in the superscript. Note that the substitution applied to the entire expression, 
$$
\left((z_1^2 - z_2^2)^r R (z_1 e^{u -v +\alpha h} / z_2)^{\pm 1}\right)\Big|_{z_1 =z_2  e^{z_0 }}\big. 
$$
does not need to be defined because for any integer $r\geqslant 0$ we have
$$
(z_1^2 - z_2^2)^r R (z_1 e^{u -v +\alpha h} / z_2)^{\pm 1}
\in
\ndo\CC^2\ot\ndo\CC^2 \ts [[h,u,v]]((z_1 /z_2 )).
$$
Hence, it is necessary to restrict  the given expression modulo $u^{a}, v^{b}, h^{k}$ first, as  in \eqref{aswith}.

Finally, we establish an analogue  of Proposition \ref{localitycor} for the series $\Uc(z)$. It can be proved by arguing as in the proof of \cite[Lemma 3.2]{K} and using \eqref{useries2}.

\begin{pro}\label{localitycor2}
\begin{enumerate}
\item
For any integer $k\geqslant 0$ there exists an integer $r\geqslant 0$ such that
$$
(1-z^2)^r \ts \Uc(z)^{\pm 1}\in
 \CC  [z^{\pm 2 } , h ]
\mod h^k .
$$
\item For any  integers $a,b,k> 0$ and $\alpha\in\CC$ there exists an integer $r\geqslant 0$ such that the coefficients of all monomials 
\beq\label{xmonomix2}
u^{a'} v^{b'} h^{k'}, \quad\text{where}\quad 0\leqslant a'\leqslant a,\, 0\leqslant b'\leqslant b,\, 0\leqslant k'\leqslant k,
\eeq
 in
$ (z_1^2 - z_2^2)^r \ts\Uc (z_1 e^{u -v +\alpha h} / z_2)^{\pm 1}$ 
belong to $ \CC [z_1^{\pm 2}   , z_2^{ \pm  2}]$ and such that the coefficients of all monomials \eqref{xmonomix2} in
$$
\left((z_1^2 - z_2^2)^r\ts \Uc (z_1 e^{u -v +\alpha h} / z_2)^{\pm 1}\right)\big|_{z_1=z_2 e^{z_0}}^{\modd u^{a}, v^{b}, h^{k}}\big. 
\fand   
z_2^{2r} (e^{2z_0}-1)^r  \ts\Uc (e^{z_0+u-v+\alpha h})^{\pm 1}
$$
coincide.

\end{enumerate}
\end{pro}

\section{Quantum vertex algebra \texorpdfstring{$\vr$}{V(R)}}\label{sec03}

In this section, we modify  the approach in \cite[Sect. 3, 4]{BK1} to generalize the construction \cite[Thm. 4.2]{BK1} of a certain quantum vertex algebra associated with the Yang $R$-matrix  to the case of    trigonometric   and  elliptic $R$-matrix $R(e^u)$.   To consider the trigonometric and elliptic setting simultaneously, we set $N = 2$ in the elliptic case; recall that we have $N\geqslant 2$ in
the trigonometric case.

Fix $c\in\CC$.
Denote by $\vr$   the topologically free  associative algebra over  $\CC[[h]]$  generated by the elements $l_{ij}^{(-r)}$, where $r=1,2\ldots  $ and $i,j=1,\ldots ,N$,   subject to the   relations 
\beq\label{rtt_trig_ell}
  R_{21}(e^{-u+v})\ts L_1^+(u)\ts L_2^+ (v)
=  \ts L_2^+ (v)\ts L_1^+(u)\ts R^*  ( e^{u-v}  ).
\eeq
The    generator matrix $L^+(u) $ is given  by
$$
L^+(u) =\sum_{i,j=1}^N e_{ij}\ot l_{ij}^+ (u),\quad\text{where}\quad l_{ij}^+(u)= \sum_{r=1}^{\infty}l_{ij}^{(-r)}u^{r-1} .
$$
Throughout the paper, we use subscripts to indicate the matrix copies in the tensor product algebras, e.g., in  $(\ndo\CC^N)^{\ot a}\ot \vr[[u]]$, we write 
$$
L_b^+(u) =\sum_{i,j=1}^N 1^{\ot (b-1)}\ot e_{ij}\ot 1^{\ot (a-b)}\ot l_{ij}^+ (u).
$$
In particular, in the defining relation \eqref{rtt_trig_ell}, we employ such a notation for
   $a=2$ and $b=1,2$. The defining relation of $\vr$ depends on the choice of $c\in\CC$ in the elliptic case, due to the term $R^*  ( e^{u-v}  )$ on the right-hand side of \eqref{rtt_trig_ell}. Nonetheless, we suppress this dependence in our notation as all results in this paper hold for any choice of $c$.

The matrix     $L^+(u)$ contains only nonnegative powers of   $u$. Furthermore,     the $R$-matrices $ R_{21}(e^{-u+v})$ and $ R^*_{12}(e^{u-v})$ possess only finitely many negative powers of $u-v $ modulo $h^n$ for all $n\geqslant 1$; recall \eqref{ptrs1} and \eqref{ptrs2}. Thus, the defining  relation \eqref{rtt_trig_ell} can be written equivalently as
\beq\label{rtt2_trig_ell}
 R_{21}(e^{v-u})\ts L_1^+(u)\ts L_2^+ (v)
=  L_2^+ (v)\ts L_1^+(u)\ts R^*_{12}(e^{-v+u}).
\eeq
Note that \eqref{rtt_trig_ell} and \eqref{rtt2_trig_ell} differ in the choice of the embedding applied on the $R$-matrices. More specifically, by the expansion convention introduced in Subsection \ref{sec0103sec0103} (recall, e.g., the map	\eqref{recall_eg}),  the identity in \eqref{rtt_trig_ell}  (resp. \eqref{rtt2_trig_ell}) employs the embedding $\iota_{u,v}$ (resp. $\iota_{v,u}$).

		From now on, we regard the matrix $L^+(u)$ as an operator series in $\ndo\CC^N \ot\vr [[u]]$, such that its action is given by the algebra multiplication.
		Denote by $\vac$   the unit in the algebra $\vr$.  
	We shall indicate by the  letter ``$h$'' in  the subscript that the corresponding $\CC[[h]]$-module is completed with respect to the $h$-adic topology. 
	For example, we write $\vr ((u))_h$ for the $h$-adic completion of $\vr((u))$, a $\CC[[h]]$-module which consists of all power series
	$a(u)=\sum_{r\in\ZZ} a_r u^{-r} \in\vr[[u^{\pm 1}]]$ such that   $\lim_{r\to \infty}a_r =0$ in the $h$-adic topology; for more information on the $h$-adic topology see \cite[Chap.  XVI]{Kas}.
	The next proposition generalizes  \cite[Prop. 3.2, 3.3]{BK1}  to the trigonometric and elliptic setting.

\begin{lem}\label{lemma53}
There exists a unique operator series
$$
L(u)\in \ndo\CC^N \ot \om(\vr,\vr((u))_h )
$$
such that for all $n\geqslant 0$ we have the identity
\begin{align}
L_0(u_0)\ts  L_1^+(u_1)\ldots L_n^+(u_n)\vac
=
 R_{0n}(e^{u_0 -u_n})^{-1}\ldots  R_{01}(e^{u_0 -u_1})^{-1}\ts
L_0^+(u_0)\ldots L_n^+(u_n)\vac\label{lemma21}
\end{align}
on  $ \ndo\CC^N \ot (\ndo\CC^N)^{\ot n}\ot \vr$ with the matrix copies labeled by $0,1,\ldots ,n$.
Moreover, the operator $L(u)$ satisfies the relation
\beq\label{rtautau}
R_{12} (e^{u-v})\ts L_1(u)\ts L_2 (v)
=
L_2 (v)\ts L_1(u)\ts R_{12}^*(e^{-v+u}).
\eeq
\end{lem}

\begin{prf}
To prove that the map $L(u_0)$ is well-defined  by \eqref{lemma21}, it suffices to check that it preserves the ideal of defining relations \eqref{rtt_trig_ell} for the algebra   $\vr$. Let $n\geqslant 2$ be an integer and $u=(u_1,\ldots ,u_n)$ a family of variables. For any $i=1,\ldots ,n-1$ consider the image of 
$ R_{i+1\ts i}(e^{-u_i +u_{i+1}})L_{1}^+(u_1)\ldots L_{n}^+(u_n)\vac$ under $L_0(u_0)$.  By \eqref{lemma21}, it  equals
\beq\label{l21}
 R_{i+1\ts i}(e^{-u_i +u_{i+1}}) \ts  R_{0n}^{-1}\ldots  R_{01}^{-1} \ts L_{0}^+(u_0)
\ts L_{1}^+(u_1)\ldots \ts L_{n}^+(u_n)\vac,
\eeq
where $R_{0i}^{-1}=R_{0i}(e^{u_0 -u_i})^{-1}$.
Due to the property \eqref{uni} (resp. 
\eqref{runi}) of the trigonometric (resp. elliptic $R$-matrix), we have
$
 R_{i+1\ts i}(e^{-u_i +u_{i+1}}) 
= R_{i\ts i+1}(e^{u_i -u_{i+1}})^{-1}\ts \Uc(e^{u_i -u_{i+1}}), 
$
so we can use the Yang--Baxter equation  \eqref{ybe2}  (resp. \eqref{rybe2}) to rewrite \eqref{l21} as
$$
    R_{0n}^{-1}\ldots  R_{0\ts i+2}^{-1} \ts  R_{0 i}^{-1}\ts  R_{0\ts i+1}^{-1}
\ts  R_{0\ts i-1}^{-1}\ldots
 R_{01}^{-1} \ts  R_{i+1\ts i}(e^{-u_i +u_{i+1}})\ts  L_{0}^+(u_0)
\ts L_{1}^+(u_1)\ldots \ts L_{n}^+(u_n)\vac.
$$
Next, by employing the defining relation \eqref{rtt_trig_ell}, this becomes
\begin{align*}
 & R_{0n}^{-1}\ldots  R_{0\ts i+2}^{-1} \ts  R_{0 i}^{-1}\ts  R_{0\ts i+1}^{-1}
\ts  R_{0\ts i-1}^{-1}\ldots
 R_{01}^{-1} \ts  \\
&\times L_{0}^+(u_0)
\ts L_{1}^+(u_1)\ldots \ts L_{i-1}^+(u_{i-1})\ts L_{i+1}^+(u_{i+1})\ts L_{i}^+(u_{i})\ts L_{i+2}^+(u_{i+2}) \ldots
\ts L_{n}^+(u_n)\vac\ts  R_{i\ts i+1}^*(e^{u_i -u_{i+1}}).
\end{align*}
Finally, we observe that this coincides with the image of
\begin{gather*}
 L_{1}^+(u_1)\ldots \ts L_{i-1}^+(u_{i-1})\ts L_{i+1}^+(u_{i+1})\ts L_{i}^+(u_{i})\ts L_{i+2}^+(u_{i+2}) \ldots
\ts L_{n}^+(u_n)\vac\ts  R_{i\ts i+1}^*(e^{u_i -u_{i+1}})\\
= R_{i+1\ts i}(e^{-u_i +u_{i+1}})L_{1}^+(u_1)\ldots L_{n}^+(u_n)\vac
\end{gather*}
under $L_0(u_0)$, thus proving the first assertion of the lemma.

The second assertion of the lemma, i.e. the identity   \eqref{rtautau}, is also verified by a straightforward computation which relies on the Yang--Baxter equation  \eqref{ybe2}  (resp. \eqref{rybe2}) and  the alternative form of the defining relation given by \eqref{rtt2_trig_ell}.
\end{prf}

	Let $n$ be a positive integer and 
   $u = (u_1 , \ldots ,u_n )$ the family of variables.
Introduce the following $R$-matrix products, with entries in
$(\ndo\CC^N )^{\ot n}  $,  
\begin{align}\label{rnm_ell_trig}
& R_{[n]}(e^{u})= \prod_{i=1,\dots,n-1}^{\longrightarrow} 
\prod_{j=i+1,\ldots,n}^{\longrightarrow} R_{ij}( e^{u_i -u_j}),
\end{align}
where the arrows indicate the order of factors. For example, if $n=3$, the expression  in \eqref{rnm_ell_trig} turns to
$R_{[3]}(e^u)=R_{12}R_{13}R_{23}$, where $R_{ij}=R_{ij}( e^{u_i -u_j})$.
Next,  for  a positive integer $m$   and $a\in\CC$ we introduce the products depending on the single variable $z$ and
the families of variables $u = (u_1 , \ldots ,u_n )$ and $v = (v_1 , \ldots ,v_m )$ with entries in  
$(\ndo\CC^N )^{\ot n} \ot (\ndo\CC^N )^{\ot m}$
 by
\begin{align}
&R_{nm}^{12}(e^{z+u-v+ah})= \prod_{i=1,\dots,n}^{\longrightarrow} 
\prod_{j=n+1,\ldots,n+m}^{\longleftarrow}   R_{ij} (e^{z+u_i-v_{j-n}+ah}).\label{oppositeof}
\end{align}
Also,  omitting the variable $z$ we have 
\beq\label{oppositeofoppositeof}
R_{nm}^{12}(e^{u-v+ah})= \prod_{i=1,\dots,n}^{\longrightarrow} 
\prod_{j=n+1,\ldots,n+m}^{\longleftarrow}   R_{ij} (e^{u_i-v_{j-n}+ah}).
\eeq
For example, if $n=m=2$, the expression in \eqref{oppositeofoppositeof} turns to
\beq\label{example_7}
R_{22}^{12}=R_{14}R_{13}R_{24}R_{23},\quad\text{where}\quad
R_{ij}=R_{ij}(e^{u_i-v_{j-n}+ah})\text{ and } R_{22}^{12}=R_{22}^{12}(e^{u-v+ah}).
\eeq
Finally, we shall write bar  in the superscript to indicate that the corresponding product of $R$-matrices comes in the order  which is opposite of \eqref{oppositeof} and \eqref{oppositeofoppositeof}, e.g.,
\begin{align*}
&R_{nm}^{\bar{\scriptstyle 1} 2}(e^{z+u-v+ah})= \prod_{i=1,\dots,n}^{\longleftarrow} 
\prod_{j=n+1,\ldots,n+m}^{\longleftarrow}   R_{ij}
, \quad
R_{nm}^{1\bar{\scriptstyle 2}}(e^{z+u-v+ah})= \prod_{i=1,\dots,n}^{\longrightarrow} 
\prod_{j=n+1,\ldots,n+m}^{\longrightarrow}   R_{ij}, \\
 &R_{nm}^{\bar{\scriptstyle 1}\bar{\scriptstyle 2}}(e^{z+u-v+ah})= \prod_{i=1,\dots,n}^{\longleftarrow} 
\prod_{j=n+1,\ldots,n+m}^{\longrightarrow}   R_{ij},
\quad\text{where}\quad 
R_{ij} =R_{ij} (e^{z+u_i-v_{j-n}+ah}).
\end{align*}
Fo example, by  using the notation from \eqref{example_7}, we get
$$R_{22}^{\bar{1}2} =R_{24}R_{23}R_{14}R_{13},\quad R_{22}^{1\bar{2}} =R_{13}R_{14}R_{23}R_{24}\fand  R_{22}^{\bar{1}\bar{2}}=R_{23}R_{24}R_{13}R_{14}.$$

 Define the following operators on $(\ndo\CC^N)^{\ot n}\ot\vr$:
$$
L^+_{[n]}(u)=L_1^+(u_1)\ldots L_n^+(u_n)
\fand 
L_{[n]}(u)=R_{[n]}(e^{u})\ts L_{1}(u_1)\ldots L_{n}(u_n).
$$
The next proposition is proved by a direct computation relying on \eqref{lemma21} and \eqref{rtautau}.
\begin{pro}
For any integers $n,m\geqslant 1$, the families of variables $u = (u_1 , \ldots ,u_n )$, $v = (v_1 , \ldots ,v_m )$ and the single variable $z$ we have
\begin{align}
&R_{nm}^{\bar{\scriptstyle 1}\bar{\scriptstyle 2}}(e^{u-v}) \ts L_{[n]}^{13}(u)\ts L_{[m]}^{23} (v)
=
L_{[m]}^{23} (v)\ts L_{[n]}^{13}(u)\ts R_{nm}^{*12}(e^{-v+u}),\non\\
&L_{[n]}^{13}(z+u)\ts L^{+23}_{[m]}(v)\vac
=
R_{nm}^{\bar{\scriptstyle 1}\bar{\scriptstyle 2}}(e^{z+u-v})^{-1}\ts 
L_{[n+m]}^+ (z+u,v)\vac,\label{lemma21gen}
\end{align}
where $z+u$  in the arguments of   $L_{[n ]}$  and $L_{[n+m]}^+  $  denotes the $n$-tuple
$(z+u_1,\ldots ,z+u_n)$.
\end{pro}

Our next goal is to equip $\vr$ with the structure of {\em quantum vertex algebra}, as defined Etingof and Kazhdan \cite[Subsect. 1.4.1]{EK}.
The following theorem extends the construction    \cite[Thm. 4.2]{BK1}   to the trigonometric and elliptic setting. It can be verified by straightforward arguments which go in parallel with its rational counterpart, so we omit the proof; see the proof of \cite[Thm. 4.2]{BK1} for details. From now on, the tensor products of $\CC[[h]]$-modules are understood as $h$-adically completed.

\begin{thm}\label{mainn}
There exists a unique   quantum vertex algebra structure on $\vr$ such that the vertex operator map   is given by
\beq\label{Ymap}
Y(L^+_{[n]}(u)\vac,z)=  L_{[n]}(z+u),
\eeq
the vacuum vector is the unit $\vac $
and the braiding  map  is defined by  
$$
\mathcal{S}(z)\big(  L_{[n]}^{+13}(u)  L_{[m]}^{+24}(v) 
 (\vac\otimes \vac) \big) 
 = R_{nm}^{\bar{\scriptstyle 1}\bar{\scriptstyle 2}}(e^{z+u-v})\ts L_{[n]}^{+13}(u) \ts  L_{[m]}^{+24}(v) \ts  R_{nm}^{*12}(e^{z+u-v})^{-1}
 (\vac\otimes \vac) .
$$
\end{thm}

\begin{rem}
As with the aforementioned result \cite[Thm. 4.2]{BK1}, Theorem \ref{mainn} closely resembles the Etingof--Kazhdan construction \cite[Thm. 2.3]{EK}; see the discussion at the end of \cite[Sect. 4]{BK1} which applies to this setting as well.  In particular, the vertex operator map \eqref{Ymap} possesses the decomposition 
$$
Y(L_{[n]}^+(u),z)=
\left( L_1^+( z+u_1)^{t_1}\ldots L_n^+( z+u_n)^{t_n}
\ts
 L_n^-( z+u_n)^{t_n}\ldots L_1^-( z+u_1)^{t_1}\right)^{t_1\ldots t_n},
$$
where $t_k$ stands for the matrix transposition $e_{ij}\to e_{ji}$ applied on the $k$-th tensor factor and $L^-(u)$ is an annihilation-like operator on $\vr$ uniquely determined by
$$
L^-_0(u_0)\ts  L_1^+(u_1)\ldots L_n^+(u_n)\vac
=
 R_{0n}(e^{u_0 -u_n})^{-1}\ldots  R_{01}(e^{u_0 -u_1})^{-1}\ts
L_1^+(u_1)\ldots L_n^+(u_n)\vac 
$$
with $n\geqslant 0$. On the other hand, the $\Sc$-locality in the Etingof--Kazhdan quantum affine VOAs is governed by the quantum current commutation relation of  Reshetikhin and Semenov-Tian-Shansky \cite{RS} (see \cite[Subsect. 2.1.4]{EK}), while in $\vr$ it takes the form of the FRT-relation \cite{FRT}.
\end{rem}

In the next corollary, we regard the   $R$-matrix as a rational function in the variable $z $; recall Remarks \ref{rationalfs_trig} and \ref{rationalfs}.  We need the map $\sigma$, established by the  corollary,  to introduce the suitable notion of     $\phi$-coordinated $\vr$-module; see    Definitions  \ref{phimod} and \ref{glavnadefinicija} below.
As with the existence of the braiding  $\Sc$ from  Theorem \ref{mainn}, the corollary can be verified by a direct argument which relies on the Yang--Baxter equations \eqref{ybe2}  and \eqref{rybe2} and the  defining relation \eqref{rtt_trig_ell}  for the algebra  $\vr$.

\begin{kor}\label{sigma_corollary} 
There exists a unique $\CC[[h]]$-module map
$$
\sigma(z)\colon \vr\ot\vr\to\vr\ot \vr \ot \CC (z)[[h]]
$$
such that we have
\beq\label{Smapratf}
\sigma(z)\big(  L_{[n]}^{+13}(u)  L_{[m]}^{+24}(v) 
 (\vac\otimes \vac) \big) 
 = R_{nm}^{\bar{\scriptstyle 1}\bar{\scriptstyle 2}}(ze^{ u-v})\ts L_{[n]}^{+13}(u) \ts  L_{[m]}^{+24}(v) \ts  R_{nm}^{*12}(ze^{ u-v})^{-1}
 (\vac\otimes \vac) .
\eeq
\end{kor}

\section{FRT-operator}\label{FRTsection} 

In this section, we introduce the notion of FRT-operator and derive its basic properties.
Its definition is motivated by the famous FRT-relation   \cite{FRT} and  also by the defining relation for the elliptic quantum algebra for $\widehat{\mathfrak{gl}}_2$ \cite{FIJKMY}. Let $c\in\CC$.

\begin{defn}\label{frtop_ell}
Let $W$ be a topologically free $\CC[[h]]$-module.
An operator series
 \beq\label{defres}
\Lc(z)=\Lc(z)_W\in\ndo\CC^N \ot \om(W, W ((z))_h)
\eeq
is said to be an {\em FRT-operator} over $W$ if it satisfies
\begin{align}
&  \iotaxy (R_{12}(x/y)) \ts \Lc_1 (x)\ts \Lc_2(y)
=
\Lc_2(y)\ts \Lc_1 (x)\ts  R^*_{12}(x/y).\label{defrelacije}
\end{align}
\end{defn}
Note that the FRT-identity  \eqref{defrelacije}  depends on the choice of $c\in\CC$ in the elliptic case because of the factor $ R^*_{12}(x/y)$; recall \eqref{zvijezdice}. 

\begin{rem}\label{rem_frt_op}
The embedding  $\iotaxy$ on the left-hand side means that the corresponding $R$-matrix $R_{12}(x/y)$ should be regarded as a rational function in the variables $x$ and $y$, as explained in Remarks \ref{rationalfs_trig} and \ref{rationalfs}, and then the embedding should be applied on the coefficients of powers of $h$. Thus, the $R$-matrix entries   produce formal power series which belong to $\CC((x))((y))[[h]]$, so that the expression on the left-hand side of \eqref{defrelacije} is well-defined.
On the other hand, the entries of the $R$-matrix $R^*_{12}(x/y)$ on the right-hand side are formal power series in $\CC((y))((x))[[h]]$. Hence, due to \eqref{defres}, the left-hand (resp. right-hand) side of \eqref{defrelacije} belongs to
$$
(\ndo\CC^N)^{\ot 2}   \ot\om(W, W ((x))((y))_h)\quad\text{(resp. }
(\ndo\CC^N)^{\ot 2}   \ot\om(W, W ((y))((x))_h)\text{ )}.
$$
Thus, the FRT-identity  \eqref{defrelacije} implies that both its sides belong to
$$
(\ndo\CC^N)^{\ot 2}   \ot\om(W, W ((x,y))_h).
$$
\end{rem}
 
Let us derive two   generalizations of  the FRT-relation \eqref{defrelacije}.
 First, in parallel with \eqref{rnm_ell_trig}, for any integer $n\geqslant 2$ and the family of variables $x=(x_1,\ldots ,x_n)$, we introduce the $R$-matrix products with entries in $(\ndo\CC^N)^{\ot n}$,
by
\begin{align*}
&R_{[n]}(x)= \prod_{i=1,\dots,n-1}^{\longrightarrow} 
\prod_{j=i+1,\ldots,n}^{\longrightarrow} R_{ij}  \Fand
 R^\prime_{[n]}(x)= \prod_{i=1,\dots,n-1}^{\longrightarrow} 
\prod_{j=i+1,\ldots,n}^{\longrightarrow}   \iotaxixj R_{ij} ,
\end{align*}
where 
$R_{ij}=R_{ij}(x_i/x_j)$  
and the arrows indicate the order of  factors.
Moreover, in the elliptic case, we write $R_{[n]}^*(x)$ and $R_{[n]}^{\prime\, *}(x)$  to indicate that the $R$-matrices $ R_{ij}(x_i/x_j)$ in the corresponding products are replaced by 
 $ R_{ij}^*(x_i/x_j)$.
The relation  \eqref{defrelacije} implies the identity of formal power series with coefficients in $(\ndo\CC^N)^{\ot n} \ot \ndo W$,
\beq\label{rttgen45}
R^\prime_{[n]}(x)\ts \Lc_1(x_1)\ldots \Lc_{n}(x_n)
=
\Lc_n(x_n)\ldots \Lc_{1}(x_1)\ts R_{[n]}^* (x).
\eeq

Denote the left-hand side of \eqref{rttgen45} by
$
\Lc_{[n]}(x)=  \Lc_{[n]}(x_1,\ldots ,x_n)$ so that we have
\beq\label{normal_trig_ell}
\Lc_{[n]}(x)=R^\prime_{[n]}(x)\ts \Lc_1(x_1)\ldots \Lc_{n}(x_n).
\eeq
 For any integers $n,m\geqslant 1$ and  the families of variables $x=(x_1,\ldots ,x_n)$ and   $y=(y_1,\ldots ,y_m)$ introduce the $R$-matrix products with entries in $(\ndo\CC^N)^{\ot n}\ot (\ndo\CC^N)^{\ot m}$,
\begin{align*}
 R_{nm}^{12}(x/y)= \prod_{i=1,\dots,n}^{\longrightarrow} 
\prod_{j=n+1,\ldots,n+m}^{\longleftarrow}  R_{ij}  \Fand
 R_{nm}^{\prime\,12}(x/y)= \prod_{i=1,\dots,n}^{\longleftarrow} 
\prod_{j=n+1,\ldots,n+m}^{\longrightarrow}\iotaxiyj R_{ij} ,
\end{align*}
where $R_{ij}=R_{ij}(x_i/y_j)$.
By using the   FRT-relation \eqref{defrelacije}, one   obtains the   identity 
of formal power series with coefficients in $(\ndo\CC^N)^{\ot n} \ot (\ndo\CC^N)^{\ot m} \ot \ndo W$, 
\beq\label{rtt7gen}
R_{nm}^{\prime\, 12}(x/y)\ts \Lc^{13}_{[n]}(x)\ts \Lc^{23}_{[m]}(y)
=
\Lc^{23}_{[m]}(y)\ts \Lc^{13}_{[n]}(x)\ts  R_{nm}^{*12}(x/y) .
\eeq
Note that the superscripts $1$ (resp. $2$ and $3$) in \eqref{rtt7gen}  indicate the tensor factors $1,\ldots ,n$ (resp. $n+1,\ldots ,n+m$ and $n+m+1$).
   Also, one can  derive analogously 
\beq\label{id886}
\Lc_{[n+m]}(x,y)=R_{nm}^{\prime \, 12}(x/y)\ts \Lc^{13}_{[n]}(x)\ts \Lc^{23}_{[m]}(y).
\eeq

By the constraint in   \eqref{defres},
all matrix entries of $\Lc(z)w$ belong to $W((z))_h$ for all $w\in W$. 
By  adapting the argument from Remark \ref{rem_frt_op}, this can be generalized as follows.

\begin{pro}\label{propres}
 Let $n\geqslant 1$ be an integer,  $x=(x_1,\ldots ,x_n)$ the family of variables and
$\Lc(z)$ an FRT-operator   over $W$. The matrix entries of
$\Lc_{[n]}(x_1,\ldots ,x_n)w$ belong to
$W((x_1 ,\ldots ,x_n ))_h$
 for all $w\in W$.
\end{pro}

\section{On \texorpdfstring{$\vr$}{V(R)}-modules in the trigonometric case}\label{sec04_trig}

In this section, we use the trigonometric setting from Subsection \ref{sec0101}. 
We use  the theory of $\phi$-coordinated modules
which was introduced by Li  \cite{Liphi}; see also \cite{JKLT,Kong} for an overview of the   theory and   more recent results. 
In general, the  notion of   $\phi$-coordinated   module can be regarded for any
associate $\phi=\phi(z_2,z_0)$ of the one-dimensional additive formal group. Throughout this paper, as in \cite[Sect. 5]{Liphi},  we    consider only the associate
$\phi(z_2,z_0) = z_2 e^{z_0}$.
In Subsection \ref{subsection_51_phi}, we show that constructing a $\phi$-coordinated $\vr$-module structure over a topologically free $\CC[[h]]$-module $W$ is equivalent to constructing an FRT-operator over $W$. Finally, in Subsection \ref{applications_52_trig}, we give some example applications of this result to certain algebras associated with the trigonometric $R$-matrix.

\subsection{\texorpdfstring{$\phi$}{phi}-coordinated modules}\label{subsection_51_phi} 

The next definition  
is based on \cite[Def. 3.4]{Liphi}, which we slightly modify  to make it compatible with Etingof--Kazhdan's quantum VOA theory;  cf. \cite[Def. 2.7, Rem. 2.8, Rem. 2.9]{K}. It employs the map $\sigma$ established by Corollary \ref{sigma_corollary}.

\begin{defn}\label{phimod}
A {\em $\phi$-coordinated $\vr$-module} is a pair $(W,Y_W)$ such that $W$ is a topologically free $\CC[[h]]$-module and $Y_W=Y_W(\cdot, z)=Y_W(z)$ is  a $\mathbb{C}[[h]]$-module map
\begin{align*}
Y_W \colon \vr\ot W&\to W[[z^{\pm 1}]]\\
u\ot w&\mapsto Y_W(z)(u\ot w)=Y_W(u,z)w=\sum_{r\in\mathbb{Z}} u_r w \ts z^{-r-1} 
\end{align*}
which satisfies the {\em vacuum property}:
 $$Y_W(\vac,z)w=w \quad\text{for all}\quad w\in W; $$
the {\em truncation condition}:
\beq\label{usual_trunc}
Y_W(v,z)\in\om(W,W((z))_h)\quad\text{for all }v\in \vr;
\eeq
the {\em weak associativity}: for any $u,v\in \vr$ and $k\in\mathbb{Z}_{\geqslant 0}$ there exists $r\in\mathbb{Z}_{\geqslant 0}$ such that
\begin{align}
&(z_1-z_2)^r\ts Y_W(u,z_1)Y_W(v,z_2)\in\om (W,W((z_1,z_2)) )\mod h^k\Fand\non\\
&\big((z_1-z_2)^r\ts Y_W(u,z_1)Y_W(v,z_2)\big)\big|_{z_1= z_2e^{z_0}}^{\modd h^k}  \big. \non\\
&\qquad- z_2^r (e^{z_0} -1)^r\ts Y_W\left(Y(u,z_0)v,z_2\right)\ts
\in\ts  h^k \om(W,W[[z_0^{\pm 1},z_2^{\pm 1}]]);\label{associativitymod_trig}
\end{align}
and
the {\em $\sigma$-locality}:
for any $u,v\in \vr$ and $k\in\mathbb{Z}_{\geqslant 0}$ there exists $r\in\mathbb{Z}_{\geqslant 0}$ such that    
\begin{align}
&(z_1-z_2)^{r}\ts Y_W(z_1)\big(1\otimes Y_W(z_2)\big)\iotaopjd\big(\sigma(z_1/z_2)(u\otimes v)\otimes w\big)\label{localitymod_trig}\\
&\qquad-(z_1-z_2)^{r}\ts Y_W(z_2)\big(1\otimes Y_W(z_1)\big)(v\otimes u\otimes w)\ts
\in\ts h^k W[[z_1^{\pm 1},z_2^{\pm 1}]] \quad\text{for all }w\in W.\non
\end{align}

Let $U$ be a topologically free $\CC[[h]]$-submodule of $W$. A pair $(U,Y_{U})$ is said to be a {\em    $\phi$-coordinated $\vr$-submodule} of $W$ if it is a    $\phi$-coordinated $\vr$-module such that $Y_U(v,z)u = Y_W(v,z)u$ for all $v\in\vr$ and $u\in U$.
\end{defn}

The next theorem is the main result of this section.

\begin{thm}\label{main_trig}
Let $\Lc(z)_W$ be an FRT-operator   over the $\CC[[h]]$-module $W$. There exists a unique structure of     $\phi$-coordinated $\vr$-module on $W$ such that
\beq\label{maineq1_trig}
Y_W(L_{[n]}^+(u_1,\ldots ,u_n)\vac , z)=\Lc_{[n]}(x_1,\ldots ,x_n)_W\big|_{x_1 =z e^{u_1 },\ldots ,\,x_n =z e^{u_n } }\quad\text{for all }n\geqslant 1.
\eeq
Conversely, let $(W,Y_W)$ be a     $\phi$-coordinated $\vr$-module.
There exists a unique FRT-operator $\Lc(z)_W$
    over $W$ such that
\beq\label{maineq2_trig}
\Lc(z)_W=Y_W(L^+(0)\vac,z  ).
\eeq
Moreover, a topologically free $\CC[[h ]]$-submodule of $W$ is a      $\phi$-coordinated $\vr$-sub\-module of $W$ if and only if it is invariant with respect to the  action \eqref{maineq2_trig} of $\Lc(z)_W$.
\end{thm}

\begin{prf}
The first assertion of the theorem can be verified by   straightforward computations which show that the expression in \eqref{maineq1_trig} defines a $\CC[[h]]$-module map which satisfies   Definition \ref{phimod}. We omit the computations as they go in parallel with the slightly more complicated elliptic case from Theorem \ref{main}, whose proof we present in  detail; see, in particular, the proofs of Lemmas \ref{lll-def}--\ref{lll-assoc} below.

Suppose $(W,Y_W)$ is a     $\phi$-coordinated $\vr$-module.
The proof of the converse relies on the Jacobi-type identity    established by Li \cite[Prop. 5.9]{Liphi},
\begin{align}
&(z_2 z)^{-1} \delta\left(\frac{z_1 -z_2}{z_2 z}\right) Y_W(z_1) (1 \ot Y_W(z_2))(u\ot v) \label{sigma_jacobi_1}\\
&\qquad -(z_2 z)^{-1} \delta\left(\frac{z_2-z_1}{-z_2 z}\right) Y_W(z_2)(1\ot Y_W(z_1))  \iotaopdj \sigma (z_2 /z_1)(v\ot u)\label{sigma_jacobi_2}\\
=&\, z_1^{-1}\delta\left(\frac{z_2 (1+z)}{z_1}\right) Y_W\left(Y(u,\log(1+z))v,z_2\right)
\qquad\text{for all }u,v\in\vr
,\label{sigma_jacobi_3}
\end{align}
where
$$
\delta(z)=\sum_{k\in\ZZ} z^k\in\CC[[z^{\pm 1}]]\fand 
\log(1+z)=-\sum_{k=1}^{\infty} \frac{(-z)^{k}}{k}\in z\CC[[z]].
$$
More precisely, even though, in contrast with \cite{Liphi}, we consider the $h$-adic setting, the identity \eqref{sigma_jacobi_1}--\eqref{sigma_jacobi_3} can be again proved by  arguing as in the proofs of \cite[Lemma 5.8]{Liphi} and \cite[Prop. 5.9]{Liphi}. Finally, one can deduce from \eqref{sigma_jacobi_1}--\eqref{sigma_jacobi_3}  the FRT-relation \eqref{defrelacije} for 
\beq\label{termss6}
\Lc(z)_W\coloneqq Y_W(L^+(0)\vac,z  )\in\ndo\CC^N\ot\om (W,W((z))_h)
\eeq
 by an argument that goes in parallel with the  proof of \cite[Lemma 3.9]{K} and relies on the properties of the $R$-matrix from Proposition \ref{prop_trig_novi}; see also \cite[Prop. 5.3]{BK1}. Nonetheless, we provide some details of the proof to cover the differences which occur in this computation due to the different form of the braiding map.

Let $n$ be a positive integer. Choose an integer $r\geqslant 0$ such that the expressions
\begin{align}
&\iotaopjd\left((z_1-z_2)^r L_{14}^+(0)\ts L_{23}^+(0)\ts R_{12}(z_2/z_1)\right)(\vac\ot\vac)\mod h^n \Fand\label{termss1}\\
&\iotaopdj\left((z_1-z_2)^r L_{14}^+(0)\ts L_{23}^+(0)\ts R_{12}(z_2/z_1)\right)(\vac\ot\vac)\mod h^n\label{termss2}
\end{align}
coincide. The application of \eqref{sigma_jacobi_1} on    \eqref{termss1} yields
\beq\label{jjac1}
(z_2 z)^{-1} \delta\left(\frac{z_1 -z_2}{z_2 z}\right)  
(z_1-z_2)^r   \Lc_2(z_1)_W\ts  \Lc_1(z_2)_W\ts \iotaopjd R_{12}(z_2/z_1).
\eeq
Next, the application of \eqref{sigma_jacobi_2} on     \eqref{termss2} yields
\beq\label{jjac2}
-(z_2 z)^{-1} \delta\left(\frac{z_2-z_1}{-z_2 z}\right)  
(z_1-z_2)^r \iotaopdj R_{12}(z_2/z_1)\ts\Lc_1(z_2)_W\ts \Lc_2(z_1)_W.
\eeq
Consider the sum of    \eqref{jjac1} and \eqref{jjac2}. Multiplying this sum by $(z_2z)^{-r}$,   using the well-known $\delta$-function identity $z^r\delta(z)=\delta(z)$
and then
 applying the residue 
$\rez_{z_2 z}$  
 we get
\beq\label{termss5}
 \Lc_2(z_1)_W\ts  \Lc_1(z_2)_W\ts \iotaopjd R_{12}(z_2/z_1)
- \iotaopdj R_{12}(z_2/z_1)\ts\Lc_1(z_2)_W\ts \Lc_2(z_1)_W.
\eeq
Finally,   applying \eqref{sigma_jacobi_3} on  \eqref{termss2} 
and   using the defining relation \eqref{rtt_trig_ell} of $\vr$, we obtain
\beq\label{termss4}
z_1^{-1}\delta\left(\frac{z_2 (1+z)}{z_1}\right) (z_1-z_2)^r\ts  Y_W (L_1^+(0)L_2^+ ( \log(1+z))\vac,z_2 ).
\eeq
Clearly, the expression $Y_W (L_1^+(0)L_2^+ ( \log(1+z))\vac,z_2 )$ contains only nonnegative powers of the variable $z$. Therefore,   multiplying   \eqref{termss4} by $(z_2z)^{-r}$,   using the  $\delta$-function identity
$$
(z_1-z_2)^r \delta\left(\frac{z_2 (1+z)}{z_1}\right) =
(z_2 z)^r
\delta\left(\frac{z_2 (1+z)}{z_1}\right) 
$$
and then
 applying the residue 
$\rez_{z_2 z}$  
 we get $0$. By combining this with \eqref{termss5}, we find
 $$ \Lc_2(z_1)_W\ts  \Lc_1(z_2)_W\ts \iotaopjd R_{12}(z_2/z_1)
= \iotaopdj R_{12}(z_2/z_1)\ts\Lc_1(z_2)_W\ts \Lc_2(z_1)_W \mod h^n .
$$
 As the integer $n$ was arbitrary, we conclude that the   equality above holds for all $n$, which means that 
\eqref{termss6} defines an FRT-operator over $W$, as required.

Finally, the  last assertion of the theorem is verified by a simple argument which relies on Proposition \ref{propres} and goes in parallel with the proof of \cite[Lemma 3.10]{K}.
\end{prf}

\begin{rem}\label{rem_54_trig}
Suppose that $\Lc(z)$ is an FRT-operator  over $W$. Then, $W$    possesses the structure of     $\phi$-coordinated $\vr$-module, as given by \eqref{maineq1_trig}. Next, using this structure, one obtains the   FRT-operator $\Lc^\prime (z)$   over $W$ via \eqref{maineq2_trig}. However, by comparing   \eqref{maineq1_trig} for $n=1$ and $u_1=0$ with \eqref{maineq2_trig}, one immediately sees that this new operator $\Lc^\prime (z)$ coincides with  $\Lc(z)$. 

Conversely, suppose that $W$ is a     $\phi$-coordinated $\vr$-module with respect to the map $Y_W$. The formula \eqref{maineq2_trig} then defines the  FRT-operator over $W$. Next, using this operator and \eqref{maineq1_trig}, one obtains a new map $Y'_W$, so that $(W,Y'_W)$ is a     $\phi$-coordinated $\vr$-module. To show that the maps $Y_W$ and $Y'_W$ coincide, we first observe that, due to the weak associativity constraint \eqref{associativitymod_trig}, $Y_W(u_r v,z)$ with $u,v\in\vr$ and $r\in\ZZ$ can be expressed in terms of $Y_W(u,z)$ and $Y_W(v,z)$. Indeed, this follows by multiplying \eqref{associativitymod_trig}  by $z_2^{-r}(e^{z_0}-1)^{-r}$. Thus, it is sufficient to check that $\vr$ is topologically spanned by the elements of any set $S$ which contains the vacuum vector $\vac$ and satisfies the requirements
\begin{align*}
 l_{ij}^{(-r)}\in S \text{ for all }  i,j\in\left\{1,\ldots ,N\right\},\,r=1,2\ldots  \fand
u_r v\in S  \text{ for all }  u,v\in S,\,r\in\ZZ.
\end{align*}
Furthermore, as $\vr$ is topologically spanned by $\vac$ and all monomials
$l_{i_1 j_1}^{(-r_1)}\ldots l_{i_m j_m}^{(-r_m)}$, it is sufficient to verify that these monomials belong to the $h$-adic completion of the $\CC[[h]]$-span of $S$. However, this can be easily proved by an inductive argument which relies on the identity \eqref{lemma21gen}. Therefore,  the maps $Y_W$ and $Y'_W$ coincide.
\end{rem}

\subsection{Applications  of Theorem \ref{main_trig}}\label{applications_52_trig} 

In this subsection, we demonstrate some applications of Theorem \ref{main_trig} to the representation theory of certain algebras associated with the $h$-Yangian $\Yg$.

\subsubsection{\texorpdfstring{$h$}{h}-Yangian extension \texorpdfstring{$\Drp$}{Dh(R)}}\label{sec0202}

	Consider the algebra $\Drp$ defined as   the $h$-adically  complete topological algebra over the ring  $\CC[[h]]$ generated by the elements $\lambda_{ij}^{\pm (r)}$, where $r=0,1,\ldots  $ and $i,j=1,\ldots ,N$, subject to the   relations
\begin{gather}
   \mathcal{L}^+_1(z_1) \ts \mathcal{L}^+_2(z_2)=
R (z_1/z_2)\ts
 \mathcal{L}^+_2(z_2)\ts  \mathcal{L}^+_1(z_1) \ts R (z_1/z_2), \label{rel1}\\
R (z_1/z_2)\ts  \mathcal{L}^-_1(z_1 ) \ts  \mathcal{L}^-_2(z_2)=
 \mathcal{L}^-_2(z_2)\ts  \mathcal{L}^-_1(z_1) \ts R (z_1/z_2),\label{rel2} \\
   \mathcal{L}^-_1(z_1) \ts  \mathcal{L}^+_2(z_2)=
	\left( R_{21} (z_2/z_1)^{t_1}\ts 
 \mathcal{L}^+_2(z_2)\ts  \mathcal{L}^-_1(z_1) ^{t_1}\right)^{t_1}.\label{rel3}
\end{gather}
The matrices of generators   $ \mathcal{L}^\pm(z)$ are given by   
\beq\label{gen_mat}
 \mathcal{L}^{\pm}(z)=\sum_{i,j=1}^N e_{ij}\ot \lambda_{ij}^\pm (z),
\eeq
where
\beq\label{gen_mat_2}
 \lambda_{ij}^+(z)=  \sum_{r=1}^{\infty}\lambda_{ij}^{+(r-1)}z^{ r-1} \fand 
\lambda_{ij}^-(z)= \delta_{ij}+h \sum_{r=0}^{\infty}\lambda_{ij}^{-(r)}z^{-r}
. 
\eeq
The algebra  $\Drp$ can be regarded as   a trigonometric counterpart of the algebra $\Dr$ introduced in
  \cite[Sect. 3]{BK1}  via defining relations governed by the Yang $R$-matrix. 

Next, define the {\em $h$-Yangian} $\Yg$ as the    associative algebra over the ring $\CC[[h]]$ generated by the elements $\lambda_{ij}^{ (r)}$, where $r=0,1,\ldots  $ and $i,j=1,\ldots ,N$,   subject to the defining  relation   \eqref{rel2}, where the    generator matrix $\mathcal{L}^-(z) $ is given by
$$
\mathcal{L}^-(z)=\sum_{i,j=1}^N e_{ij}\ot \lambda_{ij}  (z) \quad\text{with}\quad
\lambda_{ij} (z)= \delta_{ij}+h \sum_{r=0}^{\infty}\lambda_{ij}^{ (r)}z^{-r}
.
$$  
Note that the algebra $\Drp$ contains the $h$-Yangian as a subalgebra. Indeed, it is evident from the defining relations \eqref{rel1}--\eqref{rel3} that the assignments
$$
\lambda_{ij}^{+(r)}\mapsto 0, \quad  \lambda_{ij}^{-( r)}\mapsto  \lambda_{ij}^{ ( r)}, 
$$
where $i,j=1,\ldots ,N$ and $r=0,1,\ldots,$  define a homomorphism $\Drp\to\Yg$.

The following definition is motivated by the notion of restricted module  for affine Kac--Moody Lie algebra. 
\begin{defn}
A $\Drp$-module $W$ is said to be {\em restricted} if it is a topologically free $\CC[[h]]$-module and the action $\Lc^-(z)_W\in\ndo\CC^N \ot \ndo W[[z^{- 1}]]$ of the generator  matrix $\mathcal{L}^-(z)$ satisfies
$$
\mathcal{L}^-(z)_W\in\ndo\CC^N \ot \om(W,W[z^{-1}]_h)
.
$$
\end{defn}

Let us give an example of restricted $\Drp$-module. Its construction goes in parallel with \cite[Prop. 3.2]{BK1}.
Let $\ar$ be  the topologically free  associative algebra over the ring $\CC[[h]]$  generated by the elements $\lambda_{ij}^{+( r)}$, where $r=0,1, \ldots  $ and $i,j=1,\ldots ,N$,   subject to the defining  relation   \eqref{rel1}, where the    generator matrix $\mathcal{L}^+(z) $ is of the same form as in \eqref{gen_mat} and the first identity in \eqref{gen_mat_2}.

\begin{pro} \label{prop32}
There exists a unique structure of restricted $\Drp$-module over $\ar$ such that 
for all $n\geqslant 0$  the action $\mathcal{L}^\pm (u)_{\ar}$  of the generator matrices $\mathcal{L}^\pm (u)$   on
$$
 \mathcal{L}_1^+(z_1)\ldots \mathcal{L}_n^+(z_n)  \in  (\ndo\CC^N)^{\ot n}\ot \ar [[z_1,\ldots ,z_n]]
$$
is given by
\begin{align}
&\mathcal{L}^+_0(z_0  )_{\ar}\ts  \mathcal{L}_1^+(z_1)\ldots \mathcal{L}_n^+(z_n)  
=
  \mathcal{L}_0^+(z_0 )
 \mathcal{L}_1^+(z_1) \ldots \mathcal{L}_n^+(z_n),\label{asignplus} \\
&\mathcal{L}^-_0(z_0 )_{\ar} \ts  \mathcal{L}_1^+(z_1)\ldots \mathcal{L}_n^+(z_n)  
=
  R_{n0}( z_n/z_0 ) \ldots   R_{10}( z_1/z_0 ) \ts
 \mathcal{L}_1^+(z_1) \ldots \mathcal{L}_n^+(z_n). \label{asignminus}
\end{align}
\end{pro}
 
 \begin{prf}
One checks by a direct computation which relies on the defining relations \eqref{rel1}--\eqref{rel3} that the    assignments \eqref{asignplus} and \eqref{asignminus} define a structure of   $\Drp$-module over $\ar$. It is restricted due to the form of the trigonometric $R$-matrix; recall \eqref{ptrs1}.
\end{prf}

The next consequence of Theorem \ref{main_trig} extends \cite[Cor. 5.2]{BK1} to the trigonometric setting.

\begin{kor}\label{prop31}
Let $W$ be a restricted $\Drp$-module. There exists a unique structure of $\phi$-coordinated $\vr$-module over $W$   such that
\beq\label{expl_mod_map}
 Y_W(T_{[n]}^+(u)\vac,z)= 
\left( \mathcal{L}_1^+(ze^{u_1})^{t_1}_W\ldots \mathcal{L}_n^+( ze^{u_n})^{t_n}_W
\ts
 \mathcal{L}_n^-( ze^{u_n})^{t_n}_W\ldots \mathcal{L}_1^-( ze^{u_1})^{t_1}_W\right)^{t_1\ldots t_n}.
\eeq
Moreover, if $W_1\subset W$ is a  $\Drp$-submodule, then it is also a $\vr$-submodule.
\end{kor}

\begin{prf}
 To apply   Theorem \ref{main_trig}, it suffices to check that
$$
\mathcal{L}(z)_W\coloneqq   
\left( \mathcal{L}_1^+(z )^{t_1}_W  \mathcal{L}_1^-( z )^{t_1}_W\right)^{t_1 } \in\ndo\CC^N \ot\om(W,W((z))_h)
$$
satisfies the FRT-relation \eqref{defrelacije}, which can be done by arguing as in the proof of \cite[Prop. 3.3]{BK1}. Finally, the explicit expression \eqref{expl_mod_map} for the module map is an immediate consequence of \eqref{maineq1_trig} and the identity
\begin{align*}
&\Lc_{[n]}(x_1,\ldots ,x_n)_W=R^\prime_{[n]}(x_1,\ldots ,x_n)\ts \Lc_1(x_1)_W\ldots \Lc_{n}(x_n)_W\\
=&
\left( \mathcal{L}_1^+(x_1)^{t_1}_W\ldots \mathcal{L}_n^+(x_n)^{t_n}_W
\ts
 \mathcal{L}_n^-( x_n)^{t_n}_W\ldots \mathcal{L}_1^-(x_1)^{t_1}_W\right)^{t_1\ldots t_n},
\end{align*}
which follows from the defining relations \eqref{rel1}--\eqref{rel3} for  $\Drp$
and \eqref{normal_trig_ell}.
\end{prf}

\subsubsection{Generalized \texorpdfstring{$h$}{h}-Yangian \texorpdfstring{$\Yht$}{Yh(R)}}\label{sec0203}
Recall the $R$-matrices  \eqref{R2p} and \eqref{R1p}. The $R$-matrix 
$$
 \R(z)=e^{-h/2} \ts \R(z,1) = e^{-h/2} \ts g(z)^{-1}   R(z )\in\ndo\CC^N\ot\ndo\CC^N[[h]][z].
$$
 also satisfies the Yang--Baxter equation  \eqref{ybe1} and, in addition, it   possesses the property
$$
 \R_{12}(z) \ts   \R_{21}(1/z) = \left(1-e^{-h}z\right)\left(1-e^{-h}z^{-1}\right) .
$$

Consider the {\em generalized $h$-Yangian} $\Yht$,   an associative algebra  over   $\CC[[h]]$ generated by the elements 
$\tau_{ij}^{(r)}$, where $r\in\ZZ$ and $i,j=1,\ldots ,N$,
 subject to the defining relations 
\beq\label{relations}
  \R  ( z_1/z_2 )\ts  \Tc_1(z_1)\ts \Tc_2 (z_2)
=  \Tc_2 (z_2)\ts \Tc_1(z_1)\ts  \R ( z_1/z_2  ),
\eeq
where the generator matrix   $\Tc (z) $ is given by
$$
\Tc (z) =\sum_{i,j=1}^N e_{ij}\ot \tau_{ij}  (z)  \quad\text{with}\quad \tau_{ij} (z)= \sum_{r\in\ZZ} \tau_{ij}^{(-r)}z^{r }.
$$
Its definition is motivated by  \cite[Sect. 3]{FPT}. It can be regarded as a trigonometric counterpart of  the Yangian quantization of the Poisson algebra $\mathcal{O}(\mathfrak{gl}_N ((z^{-1})))$   \cite[Sect. 3]{KR}.

\begin{defn}
An $\Yht$-module $W$ is said to be {\em restricted} if it is a topologically free $\CC[[h]]$-module and the action $\Tc(z)_W\in\ndo\CC^N \ot \ndo W[[z^{\pm 1}]]$ of the generator  matrix $\Tc(z)$ is such that
$
\Tc(z)_W\in\ndo\CC^N \ot \om(W,W((z))_h)
$.
\end{defn}

The next consequence of Theorem \ref{main_trig} extends \cite[Cor. 7.1]{BK1} to the trigonometric setting.

\begin{kor}
Let $(W,Y_W)$ be a $\phi$-coordinated $\vr$-module. There exists a unique structure of restricted $\Yht$-module over $W$ such that
\beq\label{frrt}
\Tc(z)_W=Y_W(L^+(0)\vac,z).
\eeq
Moreover, if $W_1\subset W$ is a  $\vr$-submodule, then it is also a $\Yht$-submodule.
\end{kor}

\begin{prf}
By Theorem \ref{main_trig}, the expression in \eqref{frrt} defines an FRT-operator over $W$, so that $\Tc(z)_W$ satisfies \eqref{defres} and \eqref{defrelacije}. Multiplying the latter by $e^{-h/2}   g(x/y)^{-1}\in\CC[[h]][x/y]$, we see that it satisfies the defining relation \eqref{relations} for   $\Yht$. Hence, \eqref{frrt} defines a $\Yht$-module structure over $W$ which is, due to \eqref{defres}, restricted, as required.
\end{prf}

Next, we use Theorem \ref{main_trig}  to give  a  vertex algebraic interpretation of the construction \cite[Sect. 8]{BK1} of certain commutative families in the completion of $\Yht$. Let $I_p$ with $p\geqslant 0$ be the two-sided ideal in $\Yht$ generated by $h^{p+1}$ and all $\tau_{ij}^{(-r)}$ with $r>p$. Introduce the completed generalized $h$-Yangian $\Yhtc$ by
$$
\Yhtc =\lim_{\longleftarrow} \Yht / I_p.
$$
Motivated by \eqref{maineq1_trig}, we introduce the map
\begin{align}
Y_{\Yhtc} (\cdot ,z)\colon \vr\,&\mapsto\, \Yhtc[[z^{\pm 1}]]\non\\
 L_{[n]}^+(u_1,\ldots ,u_n)\vac \,&\mapsto\R_{[n]}(z e^{u_1 },\ldots ,z e^{u_n }) \ts \Tc_1(z e^{u_1 })\ldots ,\Tc_n(z e^{u_n }), \non
\end{align}
where
 $$
\R _{[n]}(x_1,\ldots ,x_n)= \prod_{i=1,\dots,n-1}^{\longrightarrow} 
\prod_{j=i+1,\ldots,n}^{\longrightarrow}   \R_{ij}(x_i/x_j) .
$$
By arguing as in the proofs of \cite[Prop. 6.1]{BK1} and \cite[Thm. 6.3]{BK1}, one can verify the following properties of the series
$$
L^+_{(n)}(u)\coloneqq \tr_{1,\ldots ,n}L_1^+(u)\ts L_2^+(u-h)\ldots L_n^+(u-(n-1)h)\vac\in\vr[[u]],\quad n=1,\ldots ,N
.$$
\begin{thm}
  The braiding map $\sigma$, as defined by   \eqref{Smapratf}, satisfies
\begin{align*}
&\sigma(z)\left( L^+_{(m)}(u)\ot L^+_{(n)}(v)\right) = L^+_{(m)}(u)\ot L^+_{(n)}(v)\quad\text{for any }m,n=1,\ldots ,N,\\
&\sigma(z)\left( a\ot L^+_{(N)}(v)\right) = a\ot L^+_{(N)}(v)\quad\text{for any }a\in\vr.
\end{align*}
\end{thm}
Finally, it remains to observe that the elements of the commutative family from \cite[Prop. 8.2]{BK1} (resp. family of central elements from \cite[Thm. 8.3]{BK1}) coincide, up to a multiplicative factor, with the coefficients of $Y_{\Yhtc} (L^+_{(n)}(0) ,z)$ with $n=1,\ldots ,N$
(resp. $Y_{\Yhtc} (L^+_{(N)}(0) ,z)$), as conjectured by \cite[Rem. 8.4]{BK1}.

\section{On \texorpdfstring{$\vr$}{V(R)}-modules in the elliptic  case}\label{sec04_ell}

In  this section, we use the elliptic setting from Subsection \ref{sec0103sec0103}.
As with   Section \ref{sec04_trig}, we  consider  the associate $ \phi(z_2,z_0)=z_2 e^{z_0}$; cf. \cite{Liphi}.
In Subsection \ref{subsection_61_ell}, we give an elliptic analogue of Theorem \ref{main_trig}, which states that constructing a generalized $\phi$-coordinated $\vr$-module structure over
a topologically free $\CC[[h]]$-module $W$ is equivalent to constructing an FRT-operator over
$W$. Finally, in Subsection \ref{sec_prf}, we prove this result.

\subsection{Generalized   \texorpdfstring{$\phi$}{phi}-coordinated modules}\label{subsection_61_ell}

The   theory of  twisted  $\phi$-coordinated modules was introduced by  Li, Tan and Wang \cite{LiTW}.
The following definition is based on the   approach and ideas from \cite{Liphi,LiTW}, which we slightly modify to adapt them to our setting.
The definition employs the map $\sigma$ established by Corollary \ref{sigma_corollary}. 

\begin{defn}\label{glavnadefinicija}
A {\em generalized    $\phi$-coordinated  $\vr$-module} is a pair $(W,Y_W)$ such that $W$ is a topologically free $\CC[[h]]$-module and $Y_W=Y_W(\cdot,z) $ is  a $\mathbb{C}[[h]]$-module map
\begin{align*}
Y_W \colon \vr\ot W&\to W[[z^{\pm 1/2}]]\\
u\ot w&\mapsto Y_W(z)(u\ot w)=Y_W(u,z)w=\sum_{r\in\frac{1}{2}\mathbb{Z}} u_r w \ts z^{-r-1} 
\end{align*}
which satisfies the {\em vacuum property}:
 \beq\label{fivac}
Y_W(\vac,z)w=w \quad\text{for all}\quad w\in W;
\eeq
the {\em $\phi$-truncation conditions}: 
\begin{align}
&Y_W(L^+(0)\vac,z )\in  \ndo\CC^2 \ot \om(W,W ((z^{1/2}))_h) \label{fitrunc2}
\end{align}
and for any positive integer $n$ and the family of variables $u=(u_1,\ldots ,u_n)$ we have
\begin{align}
&Y_W(L_{[n]}^+(u)\vac,z)\in (\ndo\CC^2)^{\ot n}\ot \om(W,W[e^{\pm u_1 },\ldots ,e^{\pm u_n }]((z^{1/2 }))_h);\label{fitrunc}
\end{align}
the {\em weak associativity}: for any $u,v\in \vr$ and $k \in\mathbb{Z}_{\geqslant 0}$ there exists $r\in\mathbb{Z}_{\geqslant 0}$ such that
\begin{align}
&(z_1-z_2)^r\ts Y_W(u,z_1)Y_W(v,z_2)\in\om (W,W((z_1^{1/2 },z_2^{1/2 })) )\mod h^k ; \label{assoc1}\\
&z_2^r (e^{2 x_0}-1)^r\ts Y_W(Y(u,x_0)v ,z_2)\in\om (W,W[e^{\pm x_0  }](( z_2^{1/2 })) )\mod h^k ; \label{assoc2}\\
&\big((z_1-z_2)^r\ts Y_W(u,z_1)Y_W(v,z_2)\big)\big|_{z_1^{1/2}=  (-1)^{t} z_2^{1/2} e^{ z_0}}^{\modd h^k}  \big. \non\\
&\qquad =
\big( z_2^r (e^{2 x_0}-1)^r\ts Y_W\left(Y(u,x_0)v,z_2\right)
\big)\big|_{x_0=  z_0+ t \pi i }^{\modd h^k }  \big.\quad\text{for }t=0,1;\label{assoc3}
\end{align}
and
the {\em $\sigma$-locality}:
for any $u,v\in \vr$ and $k \in\mathbb{Z}_{\geqslant 0}$ there exists $r\in\mathbb{Z}_{\geqslant 0}$ such that    
\begin{align}
&\quad\Big((z_1-z_2)^{r}\ts Y_W(z_1)\big(1\otimes Y_W(z_2)\big)\iotaopjdhdva\big(\sigma(z_1^{1/2}/z_2^{1/2})(u\otimes v)\otimes w\big)\Big.\label{localitymod}\\
&\Big.-(z_1-z_2)^{r}\ts Y_W(z_2)\big(1\otimes Y_W(z_1)\big)(v\otimes u\otimes w)
 \Big)\in h^k   W[[z_1^{\pm 1/2},z_2^{\pm 1/2}]]   \quad\text{for all }w\in W.\non
\end{align}

Let $U$ be a topologically free $\CC[[h]]$-submodule of $W$. A pair $(U,Y_{U})$ is said to be a {\em generalized   $\phi$-coordinated $\vr$-submodule} of $W$ if it is a generalized   $\phi$-coordinated $\vr$-module such that $Y_U(v,z)u = Y_W(v,z)u$ for all $v\in\vr$ and $u\in U$.
\end{defn}

\begin{rem}\label{remark43}
Regarding the $\phi$-truncation conditions, constraint \eqref{fitrunc2} ensures that the FRT-operator given by \eqref{maineq2} below is well-defined while \eqref{fitrunc}  implies 
\beq\label{trunc7}
Y_W(v,z)\in\om(W,W((z^{1/2}))_h)\,\text{ for all }\, v\in \vr .
\eeq
As with \eqref{usual_trunc}, the property in \eqref{trunc7} resembles the usual truncation condition, see, e.g., \cite[Def. 2.8]{LiTW}, generalized to $\CC[[h ]]$-modules, which originates from the vertex algebra theory. However, we impose  constraint \eqref{fitrunc}  instead of \eqref{trunc7} in our definition as it is required in the proof of Theorem \ref{main}; see also Remark \ref{tw_rem} below.
\end{rem}

\begin{rem}
Regarding the weak associativity property given by \eqref{assoc1}--\eqref{assoc3},  we use the ``mod'' symbol in the same way as explained in Subsection \ref{sec0103sec0103}. For example, \eqref{assoc1} means that for any $w\in W$ the expression $(z_1-z_2)^r\ts Y_W(u,z_1)Y_W(v,z_2)w$ belongs to 
$
W((z_1^{1/2 },z_2^{1/2 }))+h^k W[[z_1^{\pm 1/2 },z_2^{\pm 1/2 }]].
$
 In addition, observe that  \eqref{assoc1} and \eqref{assoc2} are necessary for the expressions in \eqref{assoc3} to be well-defined. 
Furthermore, it is worth noting that  \eqref{assoc2} and \eqref{assoc3} imply that the identity in \eqref{assoc3} holds for any integer $t$.
\end{rem}

\begin{rem}\label{tw_rem}
Unlike the notion of weak associativity which appears in the definition of {\em $\tau$-twisted $\phi$-coordinated module} \cite[Def. 2.8]{LiTW}, the constraint in \eqref{assoc3} involves $2 $ equalities for $t=0, 1$. This is due to the lack of   automorphism $\tau$ of $\vr$ such that 
$$
Y_W(\tau v ,z)=Y_W (v,z_1)\big|_{z^{1/2 }_1=- z^{1/2 }}
\quad\text{for all}\quad v\in \vr,
$$
where $(W,Y_W)$ denotes  a generalized    $\phi$-coordinated  $\vr$-module.
More specifically, in order for the ($\tau$-twisted) weak associativity  \cite[Eq. 2.19]{LiTW} to hold for $\vr$, there should exist an automorphism $\tau$    such that $\tau \colon L^+(u)\mapsto L^+(u+\pi i)$. However, the coefficients of matrix entries of $L^+(u+\pi i)$ no longer belong to $\vr$ and it is not clear how to extend this  
 quantum vertex algebra so that it accommodates such an automorphism.
Therefore, in order to avoid this problem, we   adapted  the original definition \cite[Def. 2.8]{LiTW} so that it does not involve a  quantum vertex algebra automorphism.
\end{rem}

The next theorem is the main result of this section. It is proved  by directly  verifying  the constraints from  Definitions \ref{frtop_ell} an \ref{glavnadefinicija}. 
Its proof is given in Subsection \ref{sec_prf} below.

\begin{thm}\label{main}
Let $\Lc(z)_W$ be an FRT-operator   over the $\CC[[h]]$-module $W$. There exists a unique structure of generalized   $\phi$-coordinated $\vr$-module on $W$ such that
\beq\label{maineq1}
Y_W(L_{[n]}^+(u_1,\ldots ,u_n)\vac , z)=\Lc_{[n]}(x_1,\ldots ,x_n)_W\big|_{x_1 =z^{1/2 }e^{u_1 },\ldots ,\,x_n =z^{1/2 }e^{u_n } }\quad\text{for all }n\geqslant 1.
\eeq
Conversely, let $(W,Y_W)$ be a generalized   $\phi$-coordinated $\vr$-module.
There exists a unique FRT-operator $\Lc(z)_W$
    over $W$ such that
\beq\label{maineq2}
\Lc(z)_W=Y_W(L^+(0)\vac,z^2 ).
\eeq
A topologically free $\CC[[h ]]$-submodule of $W$ is a generalized    $\phi$-coordinated $\vr$-sub\-module of $W$ if and only if it is invariant with respect to the  action \eqref{maineq2} of $\Lc(z)_W$.
\end{thm}

\begin{rem}
It is worth noting that   both conclusions of Remark \ref{rem_54_trig} can be  established  for Theorem \ref{main} by analogous arguments, which in this case rely on \eqref{assoc3}.
\end{rem}

\begin{rem}
The relation \eqref{defrelacije} in the definition of the FRT-operator can be regarded as a specialization of the defining relation for the elliptic algebra  $\aqpp$ introduced in \cite{FIJKMY}. This suggests that one can employ 
the theory of highest weight modules   developed by Foda,   Iohara,  Jimbo,   Kedem,   Miwa and Yan \cite{FIJKMY2} to construct explicit  examples of FRT-operators and, consequently,  generalized   $\phi$-coordinated $\vr$-modules.
\end{rem}

\subsection{Proof of Theorem \ref{main}}\label{sec_prf}

In this subsection, we prove Theorem \ref{main}.
Let $\Lc(z)$  be an FRT-operator   over $W$. The fact that \eqref{maineq1} defines a unique $\CC[[h]]$-module map which satisfies  weak associativity and $\sigma$-locality constraints \eqref{assoc1}--\eqref{localitymod} is proved in Lemmas \ref{lll-def}, \ref{lll-loc} and \ref{lll-assoc} below. Clearly,  the map $Y_W(\cdot, z)$ possesses the vacuum property \eqref{fivac}. Furthermore, the $\phi$-truncation conditions in \eqref{fitrunc2} and \eqref{fitrunc} follow from \eqref{defres} and Proposition \ref{propres}. Therefore, we conclude by Definition \ref{glavnadefinicija}  that \eqref{maineq1} defines a unique structure of generalized   $\phi$-coordinated $\vr $-module on $W$. The converse is proved in Lemma \ref{jac000} below.
Finally, as with the trigonometric case, the  last assertion   is verified by arguing as in the proof of \cite[Lemma 3.10]{K}
and using Proposition \ref{propres}.

\begin{lem}\label{lll-def}
The expression in  \eqref{maineq1}, together with $Y_W(\vac,z)=1_W$, defines a unique $\CC[[h]]$-module map $Y_W(\cdot ,z)$.  
\end{lem}

\begin{prf}
First, we observe that the right-hand side of \eqref{maineq1} is well-defined by Proposition \ref{propres}. Next, it is clear that \eqref{maineq1}, together with $Y_W(\vac,z)=1_W$, uniquely determines the map $Y_W(\cdot, z)$ as all coefficients of matrix entries of $L_{[n]}^+(u_1,\ldots ,u_n)\vac$ span an $h$-adically dense submodule of $\vr $.  Therefore, in order to prove the lemma, it  is sufficient to check that $Y_W(\cdot, z)$ preserves the ideal of defining relations \eqref{rtt_trig_ell}.

Let $n>i\geqslant 1$ and $k  \geqslant 0$   be arbitrary integers and $x=(x_1,\ldots ,x_n)$ the family of variables. By the first part of Proposition \ref{localitycor}, there exists an  integer $r\geqslant 0$ (which depends on $k$) such that
$$
\left(1-e^{2(u_i -u_{i+1})}\right)^r
 R_{21} (e^{-u_i+u_{i+1}}) \fand
\left(1-e^{2(u_i -u_{i+1})}\right)^r
 R_{12}^* (e^{u_i-u_{i+1}}) 
$$
possess only nonnegative powers of the variables $u_i$ and $u_{i+1}$ modulo $h^k $.
 By combining this observation with the Yang--Baxter equation \eqref{rybe1},   unitarity-like property \eqref{runi} and   the FRT-relation \eqref{defrelacije} one can verify the identity
\begin{align*}
&\left(1-e^{2(u_i -u_{i+1})}\right)^{2r}
 R_{i+1\ts i} (e^{-u_i+u_{i+1}})
\left(\Lc_{[n]}(x_1,\ldots ,x_n)\right)\big|_{x_j =z^{1/2 }e^{u_j }}\\
= &
\left(P_{i\ts i+1}\ts\Lc_{[n]}(x_1,\ldots,x_{i-1},x_{i+1},x_i,x_{i+2},\ldots ,x_n)\ts P_{i\ts i+1}\right)\big|_{x_j =z^{1/2 }e^{u_j }}\\
&\times 
\left(1-e^{2(u_i -u_{i+1})}\right)^{2r}
 R_{i\ts i+1}^* (e^{u_i-u_{i+1}}) \mod h^k,
\end{align*}
where 
the substitution $x_j =z^{1/2 }e^{u_j }$ is carried out for all $j=1,\ldots ,n$
and
$P_{i\ts i+1}$ denotes the action of the permutation operator $P\colon x\ot y\mapsto y\ot x$ applied on the tensor factors $i$ and $i+1$ of $(\ndo\CC^2)^{\ot n}$.
By the choice of the integer $r$, the given expressions contain  only nonnegative powers of the variables $u_i$ and $u_{i+1}$ modulo $h^k$. Therefore, multiplying by the inverse 
$ (1-e^{2(u_i -u_{i+1})} )^{-2r}\in\CC((u_i))[[u_{i+1}]]$ we find that
\begin{align}
& 
 R_{i+1\ts i} (e^{-u_i+u_{i+1}})
\left(\Lc_{[n]}(x_1,\ldots ,x_n)\right)\hspace{-1pt}\big|_{x_j =z^{1/2 }e^{u_j }}\Fand\label{nedj1}\\
 &
\left(P_{i\ts i+1}\ts\Lc_{[n]}(x_1,\ldots,x_{i-1},x_{i+1},x_i,x_{i+2},\ldots ,x_n)\ts P_{i\ts i+1}\right)\hspace{-1pt}\big|_{x_j =z^{1/2 }e^{u_j }}\ts
 R_{i\ts i+1}^* (e^{u_i-u_{i+1}})  \label{nedj2}
\end{align}
 coincide modulo $h^k$. As this calculation can be carried out for any choice of nonnegative integer  $k$, we conclude that the expressions in \eqref{nedj1} and \eqref{nedj2} are equal. However,  by \eqref{maineq1} we see that \eqref{nedj1} and \eqref{nedj2} coincide with the image of 
\begin{align*} 
& R_{i+1\ts i} (e^{-u_i+u_{i+1}})
\ts L_{[n]}^+ (u_1,\ldots ,u_n) 
\Fand\\
 &P_{i\ts i+1}\ts L_{[n]}^+(u_1,\ldots,u_{i-1},u_{i+1},u_i,u_{i+2},\ldots ,u_n)\ts P_{i\ts i+1} \ts
 R_{i\ts i+1}^* (e^{u_i-u_{i+1}})   
\end{align*}
under the map $v\mapsto Y_W(v,z)$. Hence, the ideal of the defining relations \eqref{rtt_trig_ell}  is preserved by $Y_W(\cdot ,z)$, i.e. the $\CC[[h]]$-module map $v\mapsto Y_W(v,z)$ is well-defined by  \eqref{maineq1}.
\end{prf}

\begin{lem}\label{lll-loc}
The map    defined $Y_W(\cdot ,z)$ possesses the
  $\sigma$-locality  property \eqref{localitymod}.
\end{lem}

\begin{prf}
Let $n$ and $m$ be nonnegative integers and $u=(u_1,\ldots ,u_n)$, $v=(v_1,\ldots ,v_m)$ the families of variables. By applying the expression 
$$Y_W(z_1)\big(1\otimes Y_W(z_2)\big)\iotaopjdhdva \sigma(z_1^{1/2}/z_2^{1/2}),$$ 
which corresponds to the first term in \eqref{localitymod}, on 
\beq\label{tn5}
L_{[n]}^{+13} (u)\ts 
L_{[m]}^{+24} (v)(\vac\ot\vac)\in
 (\ndo\CC^2)^{\ot (n+m)} \ot \vr ^{\ot 2}\ts [[u_1,\ldots ,u_n,v_1,\ldots ,v_m]]
,
\eeq
we get
\begin{align}
&\iotaopjdhdva\left( R_{nm}^{\bar{1}\bar{2}}(z_1^{1/2} e^{u-v} /z_2^{1/2})\right)
\left(\Lc_{[n]}^{13}(x_1,\ldots ,x_n)\right)\big|_{x_i =z_1^{1/2 }e^{u_i }} \non\\
 &\qquad\times\left(\Lc_{[m]}^{23}(y_1,\ldots ,y_m)\right)\big|_{y_j =z_2^{1/2 }e^{v_j }}
\ts \iotaopjdhdva\left( R_{nm}^{*12}(z_1^{1/2} e^{u-v} /z_2^{1/2})^{-1}\right),\label{tns2}
\end{align}
where, as usual, the substitutions $x_i =z_1^{1/2 }e^{u_i }$ and
$y_j =z_2^{1/2 }e^{v_j }$ are carried out for all $i=1,\ldots ,n$ and $j=1,\ldots ,m$. 
Indeed, \eqref{tns2}  follows directly from \eqref{Smapratf} and \eqref{maineq1}. Choose any  integers $a_1,\ldots ,a_n,b_1,\ldots ,b_m,k \geqslant 0$. Consider the coefficients of all monomials
\beq\label{mons1}
u_1^{a'_1}\ldots u_n^{a'_n}
v_1^{b'_1}\ldots v_m^{b'_m}
h^{k'}  , \quad \text{where}\quad
0\leqslant a'_i\leqslant a_i,\,
0\leqslant b'_j\leqslant b_j,\,
0\leqslant k'\leqslant k
\eeq
 for $i=1,\ldots ,n$, $j=1,\ldots ,m$, in \eqref{tns2}.   By    \eqref{runi}  and  the second part of Propositions \ref{localitycor} and \ref{localitycor2}, there exists an integer $r\geqslant 0$ such that the coefficients of all monomials \eqref{mons1} in
\begin{align*}
\iotaopjdhdva\left((z_1 -z_2)^r\ts  R_{nm}^{\bar{1}\bar{2}}(z_1^{1/2} e^{u-v} /z_2^{1/2})\right) 
\fand 
\iotaopdjhdva\left((z_1 -z_2)^r\ts R_{nm}^{\prime\,12}(z_1^{1/2}e^{u-v}/z_2^{1/2}  )\right)
\end{align*}
coincide. Hence the coefficients of all monomials \eqref{mons1} in the product of
$(z_1 -z_2)^r$ and \eqref{tns2} coincide with the coefficients of the corresponding monomials in
\begin{align}
&\iotaopdjhdva\left((z_1 -z_2)^r\ts  R_{nm}^{\prime\,12}(z_1^{1/2}e^{u-v}/z_2^{1/2}  )\right)
\left(\Lc_{[n]}^{13}(x_1,\ldots ,x_n)\right)\big|_{x_i =z_1^{1/2 }e^{u_i }} \label{tns3}\\
&\qquad \times  \left(\Lc_{[m]}^{23}(y_1,\ldots ,y_m)\right)\big|_{y_j =z_2^{1/2 }e^{v_j }}\ts 
\iotaopjdhdva\left( R_{nm}^{*12}(z_1^{1/2} e^{u-v} /z_2^{1/2})^{-1}\right)\label{tns4}\\
&\text{mod } u_1^{a_1},\ldots ,u_n^{a_n},v_1^{b_1},\ldots ,v_m^{b_m},h^k .\label{tns5}
\end{align}
Observe that the product of the expressions in \eqref{tns3} and \eqref{tns4} may not exist in general.
However, due to our choice of the integer $r$, all coefficients of monomials \eqref{mons1} in the given expression are well-defined as we regard the factors \eqref{tns3} and \eqref{tns4}  modulo \eqref{tns5}. Using the generalized FRT-relation \eqref{rtt7gen} we 
conclude that the  coefficients of all monomials \eqref{mons1}
in   \eqref{tns3}--\eqref{tns5} coincide with the corresponding coefficients in
\begin{align}
&\left(\Lc_{[m]}^{23}(y_1,\ldots ,y_m)\right)\big|_{y_j =z_2^{1/2 }e^{v_j }}
\left(\Lc_{[n]}^{13}(x_1,\ldots ,x_n)\right)\big|_{x_i =z_1^{1/2 }e^{u_i }} \non \\
&\qquad \times (z_1 -z_2)^r \iotaopjdhdva\left( R_{nm}^{*12}(z_1^{1/2} e^{u-v} /z_2^{1/2}) \right)
\iotaopjdhdva\left( R_{nm}^{*12}(z_1^{1/2} e^{u-v} /z_2^{1/2})^{-1}\right)\non \\
&\text{mod } u_1^{a_1},\ldots ,u_n^{a_n},v_1^{b_1},\ldots ,v_m^{b_m},h^k.\label{tns6}
\end{align}
Furthermore, by canceling the $R$-matrices we find that the coefficients of all  monomials \eqref{mons1} in the expression given by \eqref{tns6} coincide with the corresponding  coefficients in
\begin{align}
(z_1 -z_2)^r \left(\Lc_{[m]}^{23}(y_1,\ldots ,y_m )\right)\big|_{y_j  =z_2^{1/2 }e^{v_j }}
\left(\Lc_{[n]}^{13}(x_1,\ldots ,x_n )\right)\big|_{x_i =z_1^{1/2 }e^{u_i }}.\label{tns7} 
\end{align}
Finally, it remains to observe that  \eqref{tns7} equals    
$$
(z_1-z_2)^{r}\ts Y_W(z_2)\big(1\otimes Y_W(z_1)\big)
L_{[n]}^{+14} (u)\ts 
L_{[m]}^{+23} (v)(\vac\ot\vac),
$$
which corresponds to the second term in
\eqref{localitymod}, so that the $\sigma$-locality now follows.
\end{prf}

The next lemma finalizes the proof of sufficiency in Theorem \ref{main}.

\begin{lem}\label{lll-assoc}
The map   $Y_W(\cdot ,z)$ satisfies the weak associativity   \eqref{assoc1}--\eqref{assoc3}.
\end{lem}

\begin{prf}
Let $n$ and $m$ be nonnegative integers and $u=(u_1,\ldots ,u_n)$, $v=(v_1,\ldots ,v_m)$ the families of variables.
First, we verify \eqref{assoc1}. By applying $Y_W( z_1)(1\ot Y_W(z_2))$ on \eqref{tn5} and then using the relation \eqref{id886} we get
\begin{align}
&\left(\Lc_{[n]}^{13}(x )\right)\big|_{x_i =z_1^{1/2 }e^{u_i }}
\left(\Lc_{[m]}^{23}(y )\right)\big|_{y_j =z_2^{1/2 }e^{v_j }}\non\\
 = &\,
R_{nm}^{\prime\,12}(z_1^{1/2}e^{u-v}/z_2^{1/2}  )^{-1}
\left(\Lc_{[n+m]}(x,y )\right)\big|_{ x_i =z_1^{1/2 }e^{u_i },\,y_j =z_2^{1/2 }e^{v_j }}.\label{was1} 
\end{align}
 Choose any nonnegative integers $a_1,\ldots ,a_n,b_1,\ldots ,b_m,k$.
By    the second part of Propositions \ref{localitycor} and \ref{localitycor2},  there exists a nonnegative integer $r_1$ such that the coefficients of all monomials \eqref{mons1} in  the product  $(z_1-z_2)^{r_1}R_{nm}^{\prime\,12}(z_1^{1/2}e^{u-v}/z_2^{1/2}  )^{-1}$   possess only finitely many negative powers of   $z_1^{1/2 }$ and $z_2^{1/2 }$.
The same is true  for the coefficients of  \eqref{mons1} in
$$
\left(\Lc_{[n+m]}(x,y )\right)\big|_{ x_i =z_1^{1/2 }e^{u_i },\,y_j =z_2^{1/2 }e^{v_j }} \ts w
\quad\text{for any}\quad w\in W,$$
due to Proposition \ref{propres}. 
Hence, we conclude by \eqref{was1} that the constraint 
\eqref{assoc1} holds.

As for the proof of \eqref{assoc2},  the argument goes analogously. By applying $Y_W(z_2)(Y(x_0)\ot 1)$ on the expression given by \eqref{tn5} and then using \eqref{lemma21gen} and \eqref{Ymap}, we get
\begin{align}
 R_{nm}^{\bar{1}\bar{2}}(e^{x_0+u-v})^{-1}
\left(\Lc_{[n+m]}(x,y )\right)\big|_{ x_i =z_2^{1/2 }e^{ x_0+u_i },\,y_j =z_2^{1/2 }e^{v_j }}.\label{uiio}
\end{align}
Proposition \ref{propres}   implies that the coefficients of all monomials \eqref{mons1} in    
$$\left(\Lc_{[n+m]}(x,y )\right)\big|_{ x_i =z_2^{1/2 }e^{ x_0+u_i },\,y_j =z_2^{1/2 }e^{v_j }}\ts w\quad\text{for any}\quad w\in W
$$
belong to $ W[e^{\pm x_0 }]((z_2^{1/2 }))$.
Finally, due to Corollary \ref{localitycor}, there exists an integer $r>r_1$ such that the coefficients of all monomials \eqref{mons1} in
$ (e^{2x_0}-1)^r  R_{nm}^{\bar{1}\bar{2}}(e^{x_0+u-v})^{-1}$ belong to $(\ndo\CC^N)^{\ot 2}[e^{\pm x_0 }]$, which implies \eqref{assoc2}.

It remains to prove \eqref{assoc3}. 
Denote \eqref{was1} and \eqref{uiio} by $Z_1$ and $Z_2$, respectively.
By the preceding discussion, it  suffices to check that the coefficients of  \eqref{mons1} in
\begin{align}
&\left(
(z_1-z_2)^{r}\ts   Z_1
\right)\Big|^{\text{mod }u_1^{a_1},\ldots ,u_n^{a_n},v_1^{b_1},\ldots , v_m^{b_m},h^k }_{z_1^{1/2 }=(-1)^{t} z_2^{1/2 } e^{z_0  }}\qquad\text{and}\label{was2}
\\
&\left(
z_2^r(e^{2x_0}-1)^{r}\ts Z_2
\right)\Big|^{\text{mod }u_1^{a_1},\ldots ,u_n^{a_n},v_1^{b_1},\ldots , v_m^{b_m},h^k }_{x_0=z_0 +t\pi i}
\label{was3}
\end{align}
coincide for $t=0, 1$. By the choice of the integer  $r$, we can apply the above substitutions  separately on each factor and rewrite \eqref{was2} and \eqref{was3} as 
\begin{align}
&\left(
(z_1-z_2)^{r}\ts R_{nm}^{\prime\,12}(z_1^{1/2}e^{u-v}/z_2^{1/2}  )^{-1}\right)\Big|^{\text{mod }u_1^{a_1},\ldots ,u_n^{a_n},v_1^{b_1},\ldots , v_m^{b_m},h^k }_{z_1^{1/2 }=(-1)^{t} z_2^{1/2 } e^{z_0  }}\label{was4}\\
&\times\left(\left(\Lc_{[n+m]}(x,y )\right)\big|_{ x_i =z_1^{1/2 }e^{u_i },\,y_j =z_2^{1/2 }e^{v_j }}
\right)\Big|^{\text{mod }u_1^{a_1},\ldots ,u_n^{a_n},v_1^{b_1},\ldots , v_m^{b_m},h^k }_{z_1^{1/2 }=(-1)^{t} z_2^{1/2 } e^{z_0  }}\qquad\text{and} \label{was5}\\
&\left(
z_2^r(e^{2x_0}-1)^{r}\ts R_{nm}^{\bar{1}\bar{2}}(e^{x_0+u-v})^{-1}\right)\Big|^{\text{mod }u_1^{a_1},\ldots ,u_n^{a_n},v_1^{b_1},\ldots , v_m^{b_m},h^k }_{x_0=z_0 +t\pi i}\label{was6}\\
&\times\left(\Lc_{[n+m]}(x,y )\right)\Big|^{\text{mod }u_1^{a_1},\ldots ,u_n^{a_n},v_1^{b_1},\ldots , v_m^{b_m},h^k }_{ x_i =(-1)^{t} z_2^{1/2 }e^{ z_0+u_i },\,y_j =z_2^{1/2 }e^{v_j }},\label{was7}
\end{align}
respectively. Therefore, it is sufficient to prove that the coefficients of all monomials \eqref{mons1} in the given expressions \eqref{was4}--\eqref{was5} and \eqref{was6}--\eqref{was7} coincide. Clearly, \eqref{was5} is equal to \eqref{was7}. Moreover, the equality of the corresponding coefficients in \eqref{was4} and \eqref{was6} can be  established by using the property \eqref{runi}  and the second part of Propositions \ref{localitycor} and \ref{localitycor2}, so that the weak associativity constraint   \eqref{assoc3}  follows.
\end{prf}

The next lemma verifies the necessity in Theorem \ref{main}. 
Its proof employs the  identity
\begin{align}
&\delta\left(\frac{z_2(1+z)}{z_1}\right) F(z_1,z_2)  =\frac{1}{2 }\sum_{t=0,1} 
\delta\left((-1)^{t}\left(\frac{z_2(1+z)}{z_1}\right)^{\frac{1}{2} } \right)
F(z_1,z_2)\big|_{z_1^{1/2 }=(-1)^{t} (z_2(1+z))^{1/2 }},\label{dlt2}
\end{align}
which holds 
for any $F(z_1,z_2)\in\om(W,W((z_1^{1/2 },z_2^{1/2 }))_h)$, along with  the following  property of the formal delta function $\delta(z)=\sum_{r\in\ZZ}z^r$:
\begin{align}
&(z_2 z)^{-1}\delta\left(\frac{z_1 -z_2}{z_2 z}\right) -(z_2 z)^{-1}\delta\left(\frac{z_2 -z_1}{-z_2 z}\right)
=z_1^{-1}
\delta\left(\frac{z_2(1+z)}{z_1}\right);\label{dlt1}
\end{align}
see \cite[Rem. 2.8]{LiTW0} and \cite[Eq.  2.34, Eq. 2.35]{LiTW}. 

\begin{lem}\label{jac000}
Let $(W,Y_W)$ be a generalized   $\phi$-coordinated $\vr $-module.
The expression in \eqref{maineq2} defines a  FRT-operator  over $W$.
\end{lem}

\begin{prf}
Let $k $ be a nonnegative integer. By the first part of Propositions \ref{localitycor} and \ref{localitycor2}, there exists  $s\geqslant 0 $ such that  
\begin{align}
&(z_1 -z_2)^s\iotaopjdhdva\left( \Uc (z_2^{1/2}/z_1^{1/2})  R_{21}(z_2^{1/2}/z_1^{1/2})^{-1}\right) \non
\\
=\, & (z_1 -z_2)^s\iotaopdjhdva\left(\Uc (z_2^{1/2}/z_1^{1/2})  R_{21}(z_2^{1/2}/z_1^{1/2})^{-1}\right)\mod h^k,\label{prr1}
\end{align}
where, as before,  the ``mod'' symbol means that both expressions are regarded modulo $h^k$.
 Hence, by the $\sigma$-locality \eqref{localitymod}
there exists an even integer $r\geqslant 2s$ such that  
\begin{align}
&(z_1-z_2)^{r}\ts Y_W(z_1)\big(1\otimes Y_W(z_2)\big)\non
\\
& \times\iotaopjdhdva\Big( \Uc (z_2^{1/2}/z_1^{1/2})  R_{21}(z_2^{1/2}/z_1^{1/2})^{-1} L_{13}^+(0) L_{24}^+(0)(\vac\ot\vac )\ot w\Big)\non\\
=\,& (z_1-z_2)^{r}\ts Y_W(z_2)\big(1\otimes Y_W(z_1)\big) \non\\
& \times
\iotaopdjhdva
\Big(\sigma(z_2^{1/2}/z_1^{1/2})\ts\ts
\Uc (z_2^{1/2}/z_1^{1/2}) \non\\
&
\qquad\qquad\times  R_{21}(z_2^{1/2}/z_1^{1/2})^{-1} L_{14}^+(0) L_{23}^+(0)(\vac\ot\vac )\ot w\Big)\mod h^k \label{jcb1}
\end{align}
for all $w\in W$.
By the weak associativity   \eqref{assoc1}--\eqref{assoc3}, we can assume that  for all $t\in\ZZ$
\begin{align}
&(z_1 -z_2)^{r/2}\ts Y_W(L_{13}^+(0)\vac,z_1)Y_W(L_{23}^+(0)\vac,z_2)\non\\
&\qquad
\in(\ndo\CC^2)^{\ot 2}\ot\om (W,W((z_1^{1/2 },z_2^{1/2 })))\mod h^k;
\label{jcb2}\\
&z_2^{r/2} (e^{2x_0} -1)^{r/2}\ts Y_W(Y(L_{13}^+(0)\vac,x_0)L_{23}^+(0)\vac,z_2)\non\\
&\qquad\in (\ndo\CC^2)^{\ot 2}\ot\om(W,W[e^{\pm x_0 }](( z_2^{1/2 })))\mod h^k;\label{jcb2jcb2}
 \\
&\big((z_1-z_2)^{r/2}\ts Y_W(L_{13}^+(0)\vac,z_1)Y_W(L_{23}^+(0)\vac,z_2)\big)\big|_{z_1^{1/2 }=  (-1)^{t} z_2^{1/2 } e^{ z_0 }}^{\modd h^k}  \big. \non\\
&\qquad =
\big( z_2^{r/2} (e^{2 x_0}-1)^{r/2}\ts Y_W\left(Y(L_{13}^+(0)\vac,x_0)L_{23}^+(0)\vac,z_2\right)
\big)\big|_{x_0=  z_0+ t \pi i }^{\modd h^k}  \mod h^k. \big. 
\label{jcb4}
\end{align}
Indeed, if necessary, we can always choose greater even integer $r$ which, along with \eqref{jcb1}, satisfies \eqref{jcb2}--\eqref{jcb4}.
In addition, due to the second part of Proposition \ref{localitycor},   we can assume that $r$ satisfies 
\begin{align}
&\Big( (z_1-z_2)^{r/2}  \iotaopjdhdva
\Big(   R_{12}(z_1^{1/2}/z_2^{1/2})^{\pm 1}  \Big)\Big) \Big|_{z_1^{1/2 }=(-1)^{t} z_2^{1/2 } e^{z_0 }}^{\text{mod }h^k }\non  \\
 = &   \Big(
z_2^{r/2} (e^{2x_0 } -1)^{r/2}\ts
 R_{12}(e^{x_0})^{\pm 1}\Big)\Big|_{x_0 =z_0+t\pi i}^{\text{mod }h^k}
. \label{alw}
\end{align}

By applying the explicit formula  for the action of the map $\sigma$,   given by \eqref{Smapratf}, 
to    \eqref{jcb1}
and then using the property \eqref{runi}, we find that for all $w\in W$ 
 \begin{align}
&(z_1-z_2)^{r}
\iotaopjdhdva
\Big( \Uc (z_2^{1/2}/z_1^{1/2})  R_{21}(z_2^{1/2}/z_1^{1/2})^{-1}  \Big) 
Y_W(L_{13}^+(0)\vac,z_1)\ts Y_W(L_{23}^{+}(0)\vac,z_2)w \non
\\
=\,& (z_1-z_2)^{r}\ts 
 Y_W(L_{23}^{+}(0)\vac,z_2)\ts  Y_W(L_{13}^+(0)\vac,z_1)  \iotaopdjhdva \Big( R_{12}^*(z_1^{1/2}/z_2^{1/2})\Big) w  \mod h^k.\non
\end{align}
This observation implies  the   equalities, which we regard  modulo $h^k$,
\begin{align}
&(z_2 z)^{-1}\delta\left(\frac{z_1 -z_2}{z_2 z}\right) (z_2 z)^r\label{jcb9}\\
&\qquad\times 
\iotaopjdhdva
\Big( \Uc (z_2^{1/2}/z_1^{1/2}) R_{21}(z_2^{1/2}/z_1^{1/2})^{-1}  \Big) 
Y_W(L_{13}^+(0)\vac,z_1)\ts Y_W(L_{23}^{+}(0)\vac,z_2)w\non \\
&-(z_2 z)^{-1}\delta\left(\frac{z_2 -z_1}{-z_2 z}\right) (z_2 z)^r\non\\
&\qquad\times Y_W(L_{23}^{+}(0)\vac,z_2)\ts  Y_W(L_{13}^+(0)\vac,z_1)  \iotaopdjhdva \Big( R_{12}^*(z_1^{1/2}/z_2^{1/2})\Big) w\non\\
=\,&
(z_2 z)^{-1}\delta\left(\frac{z_1 -z_2}{z_2 z}\right)\bigg( (z_1-z_2)^r\non\\
&\qquad\times 
\iotaopjdhdva
\Big( \Uc (z_2^{1/2}/z_1^{1/2}) R_{21}(z_2^{1/2}/z_1^{1/2})^{-1}  \Big) 
Y_W(L_{13}^+(0)\vac,z_1)\ts Y_W(L_{23}^{+}(0)\vac,z_2)w \bigg)\non\\
&-(z_2 z)^{-1}\delta\left(\frac{z_2 -z_1}{-z_2 z}\right) \bigg((z_1-z_2)^r\non\\
&\qquad\times Y_W(L_{23}^{+}(0)\vac,z_2)\ts  Y_W(L_{13}^+(0)\vac,z_1)  \iotaopdjhdva \Big( R_{12}^*(z_1^{1/2}/z_2^{1/2})\Big) w\bigg)\non\\
=\,&
(z_2 z)^{-1}\left(
\delta\left(\frac{z_1 -z_2}{z_2 z}\right) -\delta\left(\frac{z_2 -z_1}{-z_2 z}\right)
\right)\bigg( (z_1-z_2)^r \non \\
&\qquad\times 
\iotaopjdhdva
\Big( \Uc (z_2^{1/2}/z_1^{1/2}) R_{21}(z_2^{1/2}/z_1^{1/2})^{-1}  \Big) 
Y_W(L_{13}^+(0)\vac,z_1)\ts Y_W(L_{23}^{+}(0)\vac,z_2)w \bigg).\non
\end{align}
Using the delta function identity \eqref{dlt1} we rewrite the above expression modulo $h^k $ as
\begin{align}
&z_1^{-1} 
\delta\left(\frac{z_2(1+z)}{z_1}  
\right)\bigg( (z_1-z_2)^r 
\iotaopjdhdva
\Big( \Uc (z_2^{1/2}/z_1^{1/2})  R_{21}(z_2^{1/2}/z_1^{1/2})^{-1}  \Big)\non
\\
\times & 
Y_W(L_{13}^+(0)\vac,z_1)\ts Y_W(L_{23}^{+}(0)\vac,z_2)w \bigg)\mod h^k. \label{prr2}
\end{align}

Let us rewrite \eqref{prr2}. Note that   \eqref{prr1} implies that
$$
\iotaopjdhdva
 \ts\Big( (z_1-z_2)^{r/2}\ts \Uc (z_2^{1/2}/z_1^{1/2})  R_{21}(z_2^{1/2}/z_1^{1/2})^{-1}  \Big) 
$$
belongs to $(\ndo\CC^2)^{\ot 2}((z_1^{1/2 },z_2^{1/2 }))[h]$
modulo $h^k$.
Hence, due to \eqref{jcb2}, we can apply the  identity  \eqref{dlt2} to \eqref{prr2}, thus getting
\begin{align*}
& 
\frac{1}{2 z_1}\sum_{t=0,1} 
\delta\left((-1)^{t}\left(\frac{z_2(1+z)}{z_1}\right)^{1/2 } \right)
\bigg( (z_1-z_2)^r  \iotaopjdhdva
\Big( \Uc (z_2^{1/2}/z_1^{1/2})  \\
&\qquad\times 
 R_{21}(z_2^{1/2}/z_1^{1/2})^{-1}  \Big)
Y_W(L_{13}^+(0)\vac,z_1)\ts Y_W(L_{23}^{+}(0)\vac,z_2)w \bigg)\bigg|_{z_1^{1/2 }=(-1)^{t} (z_2(1+z))^{1/2 }}^{\text{mod }h^k}.
\end{align*}
The substitution $z_1^{1/2 }=(-1)^{t} (z_2(1+z))^{1/2 }$ can be carried out in two steps,  as follows:
\begin{align}
& 
\frac{1}{2z_1}\sum_{t=0,1} 
\delta\left((-1)^{t}\left(\frac{z_2(1+z)}{z_1}\right)^{1/2 } \right) \bigg(\bigg( (z_1-z_2)^r  \iotaopjdhdva
\Big( \Uc (z_2^{1/2}/z_1^{1/2})  R_{21}(z_2^{1/2}/z_1^{1/2})^{-1}  \Big) \non\\
&\qquad\times 
Y_W(L_{13}^+(0)\vac,z_1)\ts Y_W(L_{23}^{+}(0)\vac,z_2)w \bigg)\bigg|_{z_1^{1/2 }=(-1)^{t} (z_2e^{2z_0})^{1/2 }}^{\text{mod }h^k }\bigg)\bigg|_{2z_0 =\log(1+z)}.\label{jcb6}
\end{align}

By \eqref{prr1} and \eqref{jcb2},      we have
\begin{align*}
 &\bigg( (z_1-z_2)^r  \iotaopjdhdva
\Big( \Uc (z_2^{1/2}/z_1^{1/2})  R_{21}(z_2^{1/2}/z_1^{1/2})^{-1}  \Big) \\
& \qquad\times 
Y_W(L_{13}^+(0)\vac,z_1)\ts Y_W(L_{23}^{+}(0)\vac,z_2)w \bigg)\bigg|_{z_1^{1/2 }=(-1)^{t} (z_2e^{2z_0})^{1/2 }}^{\text{mod }h^k } \\
=\,&
\bigg( (z_1-z_2)^{r/2}  \iotaopjdhdva
\Big( \Uc (z_2^{1/2}/z_1^{1/2})  R_{21}(z_2^{1/2}/z_1^{1/2})^{-1}  \Big)\bigg) \bigg|_{z_1^{1/2 }=(-1)^{t} (z_2e^{2z_0})^{1/2 }}^{\text{mod }h^k } \\
& \qquad\times \bigg(
(z_1-z_2)^{r/2}\ts
Y_W(L_{13}^+(0)\vac,z_1)\ts Y_W(L_{23}^{+}(0)\vac,z_2)w \bigg)\bigg|_{z_1^{1/2 }=(-1)^{t} (z_2e^{2z_0})^{1/2 }}^{\text{mod }h^k }. 
\end{align*}
Applying \eqref{runi} to the first and \eqref{jcb4}  to the second term on the right-hand side we get
\begin{align*}
&\bigg( (z_1-z_2)^{r/2}  \iotaopjdhdva
\Big(   R_{12}(z_1^{1/2}/z_2^{1/2})  \Big)\bigg) \bigg|_{z_1^{1/2 }=(-1)^{t} (z_2e^{2z_0})^{1/2 }}^{\text{mod }h^k }  \\
& \qquad\times \bigg(
z_2^{r/2} (e^{2x_0 } -1)^{r/2}\ts
Y_W(Y(L_{13}^+(0)\vac,x_0) L_{23}^{+}(0)\vac,z_2)w \bigg)\bigg|_{x_0 =z_0+t\pi i}^{\text{mod }h^k }. 
\end{align*}
By  the definition \eqref{Ymap} of the vertex operator map  and    \eqref{lemma21gen}, this is equal to
\begin{align*}
&\bigg( (z_1-z_2)^{r/2}  \iotaopjdhdva
\Big(   R_{12}(z_1^{1/2}/z_2^{1/2})  \Big)\bigg) \bigg|_{z_1^{1/2 }=(-1)^{t} (z_2e^{2z_0})^{1/2 }}^{\text{mod }h^k } \\
& \qquad\times \bigg(
z_2^{r/2} (e^{2x_0 } -1)^{r/2}\ts
 R_{12}(e^{x_0})^{-1}\ts
Y_W( L_{13}^+(x_0)  L_{23}^{+}(0)\vac,z_2)w \bigg)\bigg|_{x_0 =z_0 + t\pi i}^{\text{mod }h^k }. 
\end{align*}
Using $\phi$-truncation condition \eqref{fitrunc} and  \eqref{alw}  we   rewrite the given expression as
\begin{align*}
&\bigg( (z_1-z_2)^{r/2}  \iotaopjdhdva
\Big(   R_{12}(z_1^{1/2}/z_2^{1/2})  \Big)\bigg) \bigg|_{z_1^{1/2 }=(-1)^{t} (z_2e^{2z_0})^{1/2 }}^{\text{mod }h^k }  \\
& \qquad\times \bigg(
z_2^{r/2} (e^{2x_0 } -1)^{r/2}\ts
 R_{12}(e^{x_0})^{-1}\bigg)\bigg|_{x_0 =z_0+t\pi i}^{\text{mod }h^k }
\bigg(
Y_W( L_{13}^+(x_0)  L_{23}^{+}(0)\vac,z_2)w \bigg)\bigg|_{x_0 =z_0+t\pi i}^{\text{mod }h^k }  .
\end{align*}
Furthermore, due to \eqref{alw},  we can  cancel the $R$-matrices, thus getting
\begin{align*}
z_2^r (e^{2z_0} -1)^r
\Big(
Y_W( L_{13}^+(x_0)  L_{23}^{+}(0)\vac,z_2)w \Big)\Big|_{x_0 =z_0+t\pi i}^{\text{mod }h^k }. 
\end{align*}
Therefore, the expression given by \eqref{jcb6} can be written as
\begin{align}
& 
\frac{1}{2z_1}\sum_{t=0,1} 
\delta\left((-1)^{t}\left(\frac{z_2(1+z)}{z_1}\right)^{1/2 } \right)\non\\
&\qquad\times
\bigg(z_2^r (e^{2z_0} -1)^r
\Big(
Y_W( L_{13}^+(x_0)  L_{23}^{+}(0)\vac,z_2)w \Big)\Big|_{x_0 =z_0+t\pi i}^{\text{mod }h^k }\bigg)\bigg|_{2z_0 =\log(1+z)}.\label{jcb7}
\end{align}
As $\phi$-truncation condition \eqref{fitrunc} implies that the substitutions in 
\beq\label{jcb8}
\Big(\Big(
Y_W( L_{13}^+(x_0)  L_{23}^{+}(0)\vac,z_2)w \Big)\Big|_{x_0 =z_0+t\pi i}^{\text{mod }h^k}\Big)\Big|_{2z_0 =\log(1+z)}
\eeq
are well-defined, the expression given by \eqref{jcb7} can be written as
\begin{align}
& 
\frac{1}{2z_1}\sum_{t=0,1} 
\delta\left((-1)^{t}\left(\frac{z_2(1+z)}{z_1}\right)^{1/2 } \right)(z  z_2)^r \non\\
&\qquad\times
\Big(
\Big(
Y_W( L_{13}^+(x_0)  L_{23}^{+}(0)\vac,z_2)w \Big)\Big|_{x_0 =z_0+t\pi i}^{\text{mod }h^k }\Big)\Big|_{2z_0 =\log(1+z)}.\non
\end{align}

Finally,   the preceding calculation, which starts in \eqref{jcb9}, implies that the identity
\begin{align}
&(z_2 z)^{-1}\delta\left(\frac{z_1 -z_2}{z_2 z}\right) (z_2 z)^r\non\\
&\qquad\times 
\iotaopjdhdva
\Big( \Uc (z_2^{1/2}/z_1^{1/2})  R_{21}(z_2^{1/2}/z_1^{1/2})^{-1}  \Big) 
Y_W(L_{13}^+(0)\vac,z_1)\ts Y_W(L_{23}^{+}(0)\vac,z_2)w \textbf{}\non\\
&-(z_2 z)^{-1}\delta\left(\frac{z_2 -z_1}{-z_2 z}\right) (z_2 z)^r\textbf{}\non\\
&\qquad\times Y_W(L_{23}^{+}(0)\vac,z_2)\ts  Y_W(L_{13}^+(0)\vac,z_1)  \iotaopdjhdva \Big( R_{12}^*(z_1^{1/2}/z_2^{1/2})\Big) w\textbf{}\non\\
=\,&\frac{1 }{2z_1}\sum_{t=0,1} 
\delta\left((-1)^{t}\left(\frac{z_2(1+z)}{z_1}\right)^{1/2 } \right)(z  z_2)^r \non\\
&\qquad\times
\Big(
\Big(
Y_W( L_{13}^+(x_0)  L_{23}^{+}(0)\vac,z_2)w \Big)\Big|_{x_0 =z_0+t\pi i}^{\text{mod }h^k}\Big)\Big|_{2z_0 =\log(1+z)}\label{alw2}
\end{align}
holds modulo $h^k $ for all $w\in W$.
Since \eqref{jcb8} contains only nonnegative powers of  $z$,  multiplying     \eqref{alw2}  by $(zz_2)^{-r}$ and then taking the residue $\rez_{zz_2}$ we find
\begin{align}
&
\iotaopjdhdva
\Big( \Uc (z_2^{1/2}/z_1^{1/2}) R_{21}(z_2^{1/2}/z_1^{1/2})^{-1}  \Big) 
Y_W(L_{13}^+(0)\vac,z_1)\ts Y_W(L_{23}^{+}(0)\vac,z_2)w \non\\
&- Y_W(L_{23}^{+}(0)\vac,z_2)\ts  Y_W(L_{13}^+(0)\vac,z_1)  \iotaopdjhdva\Big( R_{12}^*(z_1^{1/2}/z_2^{1/2})\Big) w
=0\mod h^k .\label{jcba}
\end{align}
As the nonnegative integer $k $ and $w\in W$   were arbitrary, with $W$ being separable we conclude that the equality in \eqref{jcba} holds for all   $k\geqslant 0$, i.e. we have
\begin{align}
&
\iotaopjdhdva
\Big( \Uc (z_2^{1/2}/z_1^{1/2}) R_{21}(z_2^{1/2}/z_1^{1/2})^{-1}  \Big) 
Y_W(L_{13}^+(0)\vac,z_1)\ts Y_W(L_{23}^{+}(0)\vac,z_2)w \non\\
=\ts & Y_W(L_{23}^{+}(0)\vac,z_2)\ts  Y_W(L_{13}^+(0)\vac,z_1)  \iotaopdjhdva\Big( R_{12}^*(z_1^{1/2}/z_2^{1/2})\Big) w\quad\text{for all }w\in W.
\label{jcbajcba}
\end{align}
Finally, setting $x=z_1^{1/2}$, $y=z_2^{1/2}$   in
\eqref{jcbajcba} gives  the FRT-relation \eqref{defrelacije}, as required.
The remaining requirement \eqref{defres},
  is a consequence of the $\phi$-truncation condition \eqref{fitrunc2}.
\end{prf}

\begin{rem}
 Identity \eqref{alw2} may be regarded as a special case of the Jacobi-type identity for generalized    $\phi$-coordinated  $\vr $-modules. The Jacobi-type identity for $\tau$-twisted $\phi$-coordinated modules (recall Remark \ref{tw_rem}) was originally established in \cite[Thm. 2.13]{LiTW}. In fact, the   proof of  Lemma \ref{jac000}  is   similar  to the proof of \cite[Thm. 2.13]{LiTW}.
\end{rem}

\section{On central elements of \texorpdfstring{$\vr$}{V(R)} in the elliptic case}\label{last_section}

In this section, we consider the elliptic setting, as given in Subsection \ref{sec0103sec0103}. In Subsection \ref{subsection71}, we 
  introduce the quantum determinant for the matrix $L^+(u)$ and we show that its coefficients belong to the center of the  quantum vertex algebra $\vr$. In Subsection \ref{subsection72}, we discuss  its application  to constructing invariants of the $\phi$-coordinated   module map from Theorem \ref{main}.

\subsection{Quantum determinant} \label{subsection71}

Let $ 
A=\frac{1}{2}\left(I-P\right)
$  be the   anti-symmetrizer
on $ \CC^2 \ot \CC^2$, 
where $I$ and $P$ denote the identity and the permutation operator, respectively,  i.e.
$$I=\sum_{i,j=1,2} e_{ii}\ot e_{jj}\fand P=\sum_{i,j=1,2} e_{ij}\ot e_{ji}.$$
  Motivated by the  quantum determinant  for the type $A$ elliptic quantum algebra \cite{FIJKMY,FIJKMY2,FIR}, we define the {\em quantum determinant} of the matrix $L^+(u)$ by
\beq\label{qdet}
\qdet L^+(u) = \tr_{1, 2}\ts  A_{12} \ts L_1^+(u)\ts L_2^+(u-h/2) 
\vac,
\eeq
where the trace is taken over both copies of $\ndo\CC^2$. Clearly, 
the  quantum determinant belongs to $\vr[[u]]$.
Our goal is to show that its coefficients 
belong to the {\em center}
$$
\z(\vr)=
\left\{
v\in\vr\,:\, Y(w,z)v\in  \vr[[z]]\text{ for all }w\in \vr
\right\}
$$
of  the quantum vertex algebra $\vr$; for more information on the notion of    quantum vertex algebra center see \cite[Thm. 1.4]{DGK} and \cite[Sect. 3.2]{JKMY}.

The next lemma is well-known. Its generalization   to the {\em elliptic $R$-matrix of the $\ZZ_N$-vertex model},
 which turns to \eqref{rbar} at $N=2$,
was   established in
\cite[Lemma 5.3]{FIR}.

\begin{lem}\label{raalemma}
The following  identities hold:
\beq\label{raa}
 R_{10} (e^{u})\ts  R_{20} (e^{u-h/2}) \ts A_{12}
=   A_{12} =R_{10}^* (e^{u})\ts  R_{20}^* (e^{u-h/2}) \ts A_{12}.
\eeq
\end{lem}

\begin{prf}
 The   $R$-matrix \eqref{rbar} can be also written in terms of the Jacobi theta function; see \cite[Sect. 2.1]{FIR}. The fact that such an expression for the $R$-matrix coincides with   \eqref{rbar} is verified using the   Jacobi triple product identity. Finally, the lemma can be proved by a direct calculation which relies on   these two expressions for the $R$-matrix. 
\end{prf}

\begin{lem}
The following identities hold:
\begin{align}
&\Uc(e^u)\ts \Uc(e^{u+h/2}) 
=1, \label{uovi}\\
 &R_{02}(e^{u+ h/2})^{-1}   R_{01}(e^{u})^{-1}\ts   R_{10}(e^{-u})^{-1}
  R_{20}(e^{-u- h/2})^{-1}=1.\label{rlema}
\end{align}
\end{lem}

\begin{prf}
By using \eqref{useries2}, we rewrite the left-hand side in \eqref{uovi} as
\beq\label{uovi3}
e^{- h}  f(e^{2u})^{-1}f(e^{2u+ h})^{-1}
  f(e^{-2u })^{-1}f(e^{-2u- h})^{-1}.
\eeq
Next, by employing the identity
$f(e^u)f(e^{u+h}) =(1-e^u)(1-e^{u+ h})^{-1}$,
which follows from \cite[Eq.  2.11]{KM} for $x=e^u$ and $q=e^{h/2}$, we find that \eqref{uovi3} equals $1$, as required. 
As for the second assertion, it is sufficient to observe that, by \eqref{runi}, the left-hand side in \eqref{rlema} coincides with the inverse of the left-hand side in \eqref{uovi}.
\end{prf}

We are now ready to prove the main result of this  section.

\begin{thm}\label{qdetthm} 
The coefficients  of the quantum determinant $\qdet L^+(u)$ belong to the  center   of $\vr$. 
Moreover, we have
$$
Y(v,z)\ts \qdet L^+(u) =\qdet L^+(u)\ts Y(v,z)\quad\text{for any }v\in \vr.
$$
\end{thm}

\begin{prf}
First, we prove the equality for operators on $\ndo\CC^2\ot \vr$,
\beq\label{uovi4}
L(v_0)\ts \qdet L^+(u)
=\qdet L^+(u)\ts L(v_0),
\eeq
where $L(v_0)\in\om(\vr,\vr((v_0))_h)$ is the operator series established by Lemma \ref{lemma53}.
Let $m\geqslant 0$ be an integer and  $v=(v_1,\ldots,v_m)$ the family of variables.
Set $\wvr{u}=(u_1,u_2 )=(u,u-h/2 )$.
Applying the left-hand side of \eqref{uovi4} to $L_{[m]}^{+}(v)\vac=L_{[m]}^{+}(v_1,\ldots ,v_m)\vac$ we get
\begin{align*}
L_0 (v_0)\ts \tr_{1,2}\ts  A^1  \ts L_{[2]}^{+13}(\wvr{u})\ts
L_{[m]}^{+23}(v)
\vac
=
 \tr_{1,2}\ts  A^1 \ts L_0 (v_0)\ts L_{[2]}^{+13}(\wvr{u})\ts
L_{[m]}^{+23}(v)
\vac,
\end{align*}
where the superscripts indicate the tensor factors as follows:
$$
\smalloverbrace{\ndo\CC^2}^{0}
\ot
\smalloverbrace{(\ndo\CC^2)^{\ot 2}}^{1}
\ot
\smalloverbrace{(\ndo\CC^2)^{\ot m}}^{2}
\ot
\smalloverbrace{\vr}^{3}.
$$
Note that $L_0 (v_0)=L^{03} (v_0)$ and $A^1 =A_{12}$, so,
by \eqref{lemma21}, this equals
\begin{align*}
 \tr_{1,2}\ts  A^1\ts X\ts  R_{02}^{-1} \ts  R_{01}^{-1}\ts L_0^+ (v_0)\ts L_{[2]}^{+13}(\wvr{u})\ts
L_{[m]}^{+23}(v)
\vac, 
\end{align*}
where we use the notation $L_0^+(v_0)= L^{+03}(v_0)$,
$X=
R_{0\ts  m+2}(e^{v_0 -v_m})^{-1}\ldots R_{0\ts 3}(e^{v_0 -v_1})^{-1}$
and
$R_{0j}^{-1}=R_{0j}(e^{v_0-u_j})^{-1}$. 
Next, by using \eqref{rtt_trig_ell}, we rewrite the above expression as
\begin{align}
 \tr_{1,2}\ts  A^1\ts X\ts  R_{02}^{-1}   R_{01}^{-1}\ts R_{10}^{-1}\ts   R_{20}^{-1}\ts   L_{[2]}^{+13}(\wvr{u}) \ts  L_0^+ (v_0)
\ts  R_{02}^{*}  \ts  R_{01}^{*}
 \ts
L_{[m]}^{+23}(v)
\vac, \label{uovi5}
\end{align}
where 
$ R_{j0}^{-1}= R_{j0}(e^{-v_0+u_j})^{-1}$ and $ R_{0j}^{*}= R_{0j}^*(e^{v_0-u_j})$.
By the identity \eqref{rlema}, we have  
$
 R_{02}^{-1}  R_{01}^{-1}\ts R_{10}^{-1}\ts   R_{20}^{-1}=1,
$
so that \eqref{uovi5} equals
\begin{align}
  \tr_{1,2}\ts  A^1  \ts X\ts L_{[2]}^{+13}(\wvr{u})\ts  L_0^+ (v_0)\ts R_{02}^{*}\ts R_{01}^{*}\ts
L_{[m]}^{+23}(v)
\vac.  \label{uovi6}
\end{align}
However, by combining \eqref{runi} and \eqref{uovi} we get
$$
R_{02}^{*}\ts R_{01}^{*}
=
(R_{20}^{*})^{-1} \ts (R_{10}^{*})^{-1}\ts
\Uc (e^{v_0-u_2})\ts \Uc (e^{v_0-u_1})
= (R_{20}^{*})^{-1}\ts (R_{10}^{*})^{-1}
$$
for $(R_{j0}^{*})^{-1}=R_{j0}^{*}(e^{-v_0 +u_j})^{-1}$. 
Therefore,  \eqref{uovi6} coincides with
\beq\label{kndl}
  \tr_{1,2}\ts  A^1  \ts X\ts L_{[2]}^{+13}(\wvr{u})\ts  L_0^+ (v_0)\ts (R_{20}^{*})^{-1}\ts (R_{10}^{*})^{-1}\ts
L_{[m]}^{+23}(v)
\vac.   
\eeq
Note that we have
$
A^1\ts
L_{[m]}^{+23}(v)
=
L_{[m]}^{+23}(v)\ts A^1
$ and $
X\ts L_{[N]}^{+13}(\wvr{u})
=
L_{[N]}^{+13}(\wvr{u})\ts X
$,
as  these operators act on different tensor factors.
Hence, using the cyclic property of the trace to move the anti-symmetrizer in \eqref{kndl} to the right
and then employing \eqref{raa} we find 
\begin{align*}
   &  \tr_{1,2}  \ts X\ts L_{[2]}^{+13}(\wvr{u})\ts  L_0^+ (v_0)\ts (R_{20}^{*})^{-1}\ts (R_{10}^{*})^{-1}\ts  A^1\ts
L_{[m]}^{+23}(v)
\vac    \\
=&\,
\tr_{1,2}  \ts X\ts L_{[2]}^{+13}(\wvr{u})\ts  L_0^+ (v_0)\ts    A^1\ts
L_{[m]}^{+23}(v)
\vac   
= 
\tr_{1,2}  \ts L_{[2]}^{+13}(\wvr{u})\ts X\ts  L_0^+ (v_0)\ts
L_{[m]}^{+23}(v)\ts    A^1
\vac.
\end{align*}
Finally, by  \eqref{lemma21}, we have
$X  L_0^+  (v_0)  L_{[m]}^{+23}(v)\vac=L_0  (v_0)  L_{[m]}^{+23}(v)\vac  $,
so, returning the anti-symmetrizer to the left, we obtain
$$
\tr_{1,2}  \ts    A^1\ts L_{[2]}^{+13}(\wvr{u})\ts L_0  (v_0)\ts
L_{[m]}^{+23}(v)
\vac
=\qdet L^+(u) \ts L_0  (v_0)\ts
L_{[m]}^{+23}(v)\vac.
$$
This coincides with the action of the
  right-hand side of \eqref{uovi4}   on $L_{[m]}^{+}(v)\vac$, so we conclude that the equality in \eqref{uovi4} holds.

Finally, consider the expression
\beq\label{uovi9}
Y(L^+_{[m]}(v)\vac,z) \ts \qdet L^+(u),\quad\text{where}\quad m\geqslant 0\fand v=(v_1,\ldots ,v_m).
\eeq
as an operator  on $\vr$.
By combining  \eqref{Ymap} 
and 
\eqref{uovi4},
we conclude that
\eqref{uovi9}   equals
\beq\label{uovia}
 L_{[m]}(z+v)\ts \qdet L^+(u)= \qdet L^+(u) \ts L_{[m]}(z+v)= \qdet L^+(u) \ts Y(L^+_{[m]}(v)\vac,z).
\eeq
As $Y(w,z)\vac$ belongs to $\vr[[z]]$ for any $w\in \vr$,
the term $Y(L^+_{[m]}(v)\vac,z)\vac$ does not contain any negative powers of the variable $z$. Hence, by applying the expression in \eqref{uovi9}, which equals  \eqref{uovia}, to the vacuum vector $\vac$, we conclude that the coefficients of the quantum determinant belong to the  center  $\z(\vr)$.
Moreover,  the equalities of operators on $\vr$ in \eqref{uovi9} and \eqref{uovia} imply the second assertion of the theorem.
\end{prf}

\subsection{Submodule of invariants}\label{subsection72}

Let $\mathcal{L}(z)=\mathcal{L}(z)_W$ be an FRT-operator over $W$ and $Y_W(\cdot ,z)$ 
the corresponding generalized $\phi$-coordinated $\vr$-module map defined as in Theorem \ref{main}. We will show that the  quantum determinant \eqref{qdet} can be employed to  construct      elements of the {\em submodule of invariants} of $W$, 
\beq\label{invdef}
\z(W)=
\left\{
w\in W\ts:\ts \Lc_{[n]}(z_1,\ldots ,z_n)  w\in (\ndo\CC^2)^{\ot n} \ot W[[z_1,\ldots ,z_n]]\text{ for }n\geqslant 1
\right\}.
\eeq

\begin{thm}\label{invsthm}
For any   $w\in \z(W)$ all coefficients of the series 
\begin{align}
&Y_{W}(\qdet L^+(0)\vac , z^2)w =
\tr_{1,2} \ts A_{12} \ts \Lc_{[2]}(x_1, x_2)w  \big|_{x_1 =z ,\ts x_2 =z e^{-h/2 } }\label{rhskor} 
\end{align}
 belong
to the  submodule of invariants
$\z(W)$. 
\end{thm}

\begin{prf}
In order to simplify the notation,
we introduce the families of variables
$ 
x=(x_1,x_2 )$  and $\bar{z}=(z,ze^{-h/2} )
$ 
and denote the right-hand side of \eqref{rhskor} by
\beq\label{rhskor2}
\tr_{1,2}\ts  A^1\ts \Lc_{[2]}(x)w\big|_{x=\bar{z}} 
= \tr_{1,2}\ts  A^1\ts \Lc_{[2]}(\bar{z})w.
\eeq
By \eqref{invdef}, it is sufficient to check that, by applying $\Lc_{[n]}(y)=\Lc_{[n]}(y_1,\ldots ,y_n)$ on \eqref{rhskor2}, we obtain only nonnegative powers of the variables $y_1,\ldots ,y_n$. First, we have
$$
\Lc_{[n]}^{23}(y)\ts \tr_{1,2}\ts  A^1\ts \Lc_{[2]}^{13}(\bar{z})w
=\tr_{1,2}\ts  A^1\ts \Lc_{[n]}^{23}(y)\ts  \Lc_{[2]}^{13}(\bar{z})w,
$$
where the trace, the anti-symmetrizer $A=A_{12}=A^1$ and $\Lc_{[2]}(\bar{z})=\Lc_{[2]}^{13}(\bar{z})$ are applied on the tensor factors $1,2$ and $\Lc_{[n]}(y)=\Lc_{[n]}^{23}(y)$ is applied on the  tensor  factors $3,\ldots ,2+n$ of
$ (\ndo\CC^2)^{\ot 2}\ot (\ndo\CC^2)^{\ot n}\ot W.$  
By combining \eqref{rtt7gen} and \eqref{id886}, we get
$$
\tr_{1,2}\ts  A^1\ts \Lc_{[n]}^{23}(y)\ts  \Lc_{[2]}^{13}(\bar{z})w
=
\tr_{1,2}\ts  A^1\ts  \Lc_{[2+n]}(\bar{z},y)w\ts
 R_{2\ts n}^{*12}(\bar{z}/y)^{-1}.
$$
Next, we use the cyclic property of the trace to move the anti-symmetrizer   to the right:
\beq\label{rhskor3}
\tr_{1,2}\ts      \Lc_{[2+n]}(\bar{z},y)w\ts
 R_{2\ts n}^{*12}(\bar{z}/y)^{-1}  A^1
=
\tr_{1,2}\ts   \Lc_{[2+n]}(\bar{z},y)w\ts
   A^1.
\eeq
Note that the equality in \eqref{rhskor3} follows from the identity
$R_{2\ts n}^{*12}(\bar{z}/y)^{-1}  A^1=A^1$, which is a direct consequence of
\eqref{raa}. 
Finally,  by  \eqref{invdef},
the expression
$$
\Lc_{[2+n]}( \bar{z},y)w=
\Lc_{[2+n]}( x_1,\ldots ,x_2,y_1,\ldots ,y_n)w\big|_{x_1 =z ,\ts x_2 =z e^{-h/2} },
$$
and, consequently, the entire right-hand side of \eqref{rhskor3},
possesses only nonnegative powers of the variables $y_1,\ldots ,y_n$. Hence,  the coefficients of \eqref{rhskor}
belong to the submodule of invariants, as required.
\end{prf}

\begin{rem}
The form of series \eqref{rhskor}   resembles the action of the quantum determinant for the   elliptic algebra (cf. \cite{FIJKMY,FIJKMY2,AFRS,FIR}), which
 indicates a possible interpretation of the family of central elements for the elliptic   algebra at the critical level, as given by  Avan, Frappat, Rossi and  Sorba  \cite[Thm. 2]{AFRS}, in terms of quantum vertex algebra theory.
\end{rem}

\begin{rem}
The results from Subsections 
\ref{subsection71} and \ref{subsection72} do not immediately extend to the trigonometric setting. Namely, the trigonometric $R$-matrix $R(e^u)$, as defined by \eqref{R1p}, does not seem to  possess  the property \eqref{raa}, which is essential in the proof  of Theorem \ref{qdetthm}. However, in parallel with Section \ref{sec03}, one can construct the  quantum vertex algebra associated with a suitable renormalization of $R(e^u)$ which satisfies the trigonometric analogue of the  identity \eqref{raa}.
This  renormalization can be obtained by the use of the power series \eqref{in3}. However, in contrast with the elliptic setting, the trigonometric analogue of \eqref{raa} will accommodate the action of the   anti-symmetrizer coming from the action of the symmetric group $\mathfrak{S}_n$ of order $n$ on   $(\CC^N)^{\ot n}$, determined by the requirement that the transpositions $(i,i+1)\in\mathfrak{S}_n $ act  as the {\em $h$-permutation operator}  
$$
P^h = \sum_{i=1}^N e_{ii}\ot e_{ii} + e^{h/2}\sum_{\substack{i,j=1\\i> j}}^N e_{ij}\ot e_{ji} +e^{-h/2}\sum_{\substack{i,j=1\\i< j}}^N e_{ij}\ot e_{ji}  
$$
on the tensor factors $i$ and $i+1$, i.e. such that we have $(i,i+1)=P^h_{i\ts i+1}$.
\end{rem}

\section*{Acknowledgement} 
L.B. is member of Gruppo Nazionale per le Strutture Algebriche, Geometriche e le loro Applicazioni  (GNSAGA) of the Istituto Nazionale di Alta Matematica (INdAM).
This work has been supported  by Croatian Science Foundation under the project UIP-2019-04-8488. Furthermore, this work was supported by the project ``Implementation of cutting-edge research and its application as part of the Scientific Center of Excellence for Quantum and Complex Systems, and Representations of Lie Algebras'', PK.1.1.02, European Union, European Regional Development Fund.

\linespread{1.0}

\end{document}